\documentclass[final]{siamltex}
\usepackage{algorithm,algpseudocode}
\usepackage{latexsym, graphicx, epsfig, amsmath, amsfonts,amssymb,bm}
\usepackage{caption, subcaption}
\usepackage{epstopdf}
\usepackage{booktabs}
\usepackage{array}
\usepackage{color}
\usepackage{lineno}
\usepackage{url}

\def\be{\begin{equation}}
\def\ee{\end{equation}}

\def\x{\mathbf{x}}
\def\hx{\widehat{\x}}
\def\xta{\mathbf{x}(t; \aalpha)}
\def\xty{\mathbf{x}(t; \y)}
\def\hxta{\wh{\mathbf{x}}(t; \aalpha)}
\def\hxty{\wh{\mathbf{x}}(t; \y)}

\def\y{\mathbf{y}}
\def\f{\mathbf{f}}

\def\z{\mathbf{z}}

\def\u{\mathbf{u}}
\def\v{\mathbf{v}}
\def\aalpha{\boldsymbol{\alpha}}
\def\PPhi{\boldsymbol{\Phi}}
\def\PPsi{\boldsymbol{\Psi}}

\def\N{\mathbf{N}}
\def\hN{\widehat{\mathbf{N}}}
\def\W{\mathbf{W}}

\def\r{\rho_{\aalpha}}
\def\tr{\widetilde{\rho}_{\aalpha}}
\def\tC{\widetilde{C}}

\def\E{{\mathbb E}}
\def\Ea{\mathbb{E}_{\aalpha}}
\def\tEa{\widetilde{\mathbb{E}}_{\aalpha}}
\def\Eax{{\mathbb E}_{\aalpha}[\x(t; \aalpha)]}
\def\Eahx{{\mathbb E}_{\aalpha}[\wh{\x}(t; \aalpha)]}
\def\tEahx{\widetilde{\mathbb E}_{\aalpha}[\wh{\x}(t; \aalpha)]}
\def\Var{\mathrm{Var}}

\def\Vax{\mathrm{Var}_{\aalpha}[\x(t; \aalpha)]}
\def\Vahx{\mathrm{Var}_{\aalpha}[\wh{\x}(t; \aalpha)]}
\def\tVahx{\widetilde{\mathrm{Var}}_{\aalpha}[\wh{\x}(t; \aalpha)]}
\def\R{{\mathbb R}}

\def\Ia{I_{\aalpha}}

\def\tIa{\widetilde{I}_{\aalpha}}

\newcommand{\argmin}{\operatornamewithlimits{argmin}}

\newcommand{\wh}[1]{\widehat{#1}}

\newcommand{\figref}[1]{Fig.~\ref{#1}}
\newcommand{\norm}[1]{\left\lVert#1\right\rVert}
\newcommand{\abs}[1]{\left\lvert#1\right\rvert}

\title{Deep learning of parameterized equations with applications to uncertainty quantification}

\author{Tong Qin\footnotemark[1] \and Zhen Chen\footnotemark[1] 
	\and John D. Jakeman\thanks{Optimization and Uncertainty Quantification Department, Sandia National Laboratory, Albuqerque, NM, 87123 USA ({\tt Email: jdjakeman@sandia.gov}).} 
	\and Dongbin Xiu\thanks{Department of Mathematics,
			The Ohio State University, Columbus, OH 43210, USA.
			({\tt Emails: qin.428@osu.edu, 
				chen.7168@osu.edu, 
				xiu.16@osu.edu}).
			Funding: This work was partially supported by AFOSR FA9550-18-1-0102. }
}

\begin{document}
\maketitle
\begin{abstract}
We propose a numerical method for discovering unknown parameterized dynamical systems by using observational data of the state variables. Our method is built upon and extends the recent work of discovering unknown dynamical systems, in particular those using deep neural network (DNN).  We propose a DNN structure, largely based upon the residual network (ResNet), to not only learn the unknown form of the governing equation but also take into account the random effect embedded in the system, which is generated by the random parameters. Once the DNN model is successfully
constructed, it is able to produce system prediction over longer term and for arbitrary parameter values. For uncertainty quantification, it allows us to conduct uncertainty analysis by evaluating solution statistics over the parameter space. 
\end{abstract}

\section{Introduction} \label{sec:intro}

The ability to construct predictive models from data is essential to most if not all quantitative scientific and engineering disciplines.   
A significant amount of research has been conducted in this direction. Early efforts include symbolic regression
(\cite{bongard2007automated,schmidt2009distilling}), equation-free
modeling \cite{kevrekidis2003equation}, heterogeneous multi-scale
method (HMM) (\cite{E_HMM03}), 
artificial neural networks (\cite{gonzalez1998identification}),
nonlinear regression (\cite{voss1999amplitude}), 
empirical dynamic modeling (\cite{sugihara2012detecting,ye2015equation}), 
nonlinear Laplacian spectral analysis (\cite{giannakis2012nonlinear}), 
automated inference of dynamics
(\cite{schmidt2011automated,daniels2015automated,daniels2015efficient}),
etc. More recent effort cast the problem as one of function approximation problems and employ data-driven approaches to resolve it. The unknown governing equations are treated as target functions, which map state variables to their time derivatives. One popular approach is to employ sparsity-promoting algorithms, such as the least absolute shrinkage and selection operator (LASSO) \cite{tibshirani1996regression}, to select basis functions from a large dictionary set, which contains all candidate models, see for example \cite{brunton2016discovering, tran2017exact}. Numerous other approaches have also been developed such as those based upon projection into polynomial spaces \cite{WuXiu_JCPEQ18, wu2019structure}, dynamic mode decomposition \cite{schmid2010dynamic}, operator inference \cite{peherstorfer2016data}, model selection approach  \cite{Mangan20170009}, and 
Gaussian process regression \cite{raissi2017machine}.

Recent progress in machine learning, especially in deep neural networks (DNN), has provided new perspectives to data-driven modeling. Approaches such as physics informed neural networks \cite{raissi2017machine, raissi2018deep}, PDE-net \cite{long2018pde, long2019pde}, DNN with Runge-Kutta or multi-step integrator structures \cite{raissi2018multistep, rudy2019deep}, etc. aim at recovering hidden models in differential equation forms.
Different network structures and properties are explored in the context of recovering ODE systems \cite{rudy2019deep, raissi2018multistep, QinWuXiu2019} and PDEs  (\cite{long2017pde, long2018pde,WuXiu_JCP20}). 
DNNs have been used in other aspects of 
scientific computing. These include, construction of reduced order models \cite{HesthavenU_JCP18}, aiding numerical solvers of conservation laws (\cite{RayHeasthaven_JCP18, ray2019}), approximation of Koopman operator (\cite{brunton2017chaos}), solving differential equations \cite{lagaris1998, han2018, weinan2018deep, sirignano2018dgm}, etc. 
For uncertainty quantification (UQ).
DNNs were employed to approximate mappings from random parameters to quantities of interest (QoI) in \cite{tripathy2018}; incorporated in Bayesian framework and used as surrogate models \cite{zhu2018, zhu2019}; adopted to approximate distribution of  QoI \cite{yang2019}, etc.

The focus of this paper is on recovery of unknown parameterized dynamical systems. That is, not only is the form of the governing equations unknown, the system also possesses parameters that are unknown. Systems with unknown/uncertain parameters are often studied in the context of uncertainty quantification (UQ). When the governing equations are known, various standard UQ techniques, e.g, those based on generalized polynomial chaos \cite{ghanem1990,xiu2002}, can be readily applied.
These techniques are obviously not applicable when the governing equations are unavailable. In this paper, we poposed a data driven method that uses observation
data of the state variables to recover unknown governing equations with embedded unknown/uncertain parameters. Our method is based on and extends the work of
\cite{QinWuXiu2019}. In \cite{QinWuXiu2019}, residual network (ResNet) (\cite{he2016deep}) was used to recover unknown {\em deterministic} dynamical systems.
This paper extends the approach to unknown stochastic/random dynamical systems. To accomplish this, we introduce additional inputs in the DNN structure to
incorporate the unknown system parameters, as well as an input for ``time step''.  The proposed DNN is thus able to register system
responses with respect to different system parameters. This in turn allows us to create DNN model for the evolution of the underlying unknown equations.
The introduction of the time-step input allows the network to incorporate measurement data over non-uniform time levels. The adds more flexibility in the DNN
modeling construction. It also allows us to conduct system prediction using different time steps.
Note that our DNN modeling utilizes approximation of the flow map of the unknown system. This was proposed in \cite{QinWuXiu2019} and is different from many
other recovery methods. 
This flow-map approach has two advantages: (i) the approximation is based exact time integration and induces no temporal error associated with time step; (ii) temporal numerical derivatives of the observational data, which are usually sensitive to noises, are not required. 

Once the DNN models are successfully constructed, they allow us to conduct system predictions over longer term. More importantly, the proposed method is able to
explore the parameter space and create model prediction at arbitrary parameter values. Uncertainty quantification of the unknown system thus becomes a post-process, as we can
conduct UQ on the recovered DNN model using proper sampling methods. To this end, we also conduct some theoretical analysis to quantify the numerical errors in such UQ prediction.

%

\section{Setup and Preliminaries} \label{sec:setup}

Let us consider a parameterized system
\begin{equation}
\label{eq:ode}
\frac{d}{dt} \x (t; \aalpha)=\f(\x, \aalpha), \qquad \x(0)=\x_0,
\end{equation}
where $\x=(x_1, \dots, x_d)\in I_\x\subseteq \R^d$ are state variables and $\aalpha=(\alpha_1, \dots, \alpha_\ell)\in I_{\aalpha}\subseteq \R^\ell$ are system parameters. We are interested in the solution behavior with respect to varying parameters. In the context of uncertainty quantification (UQ), which is a major
focus of this paper, the parameters are equipped with a probability measure over $ I_{\aalpha}$. We are interested in understanding the various
solution statistics with respect to the input ${\aalpha}$.

The basic assumption of this paper is that the form of the governing equations \eqref{eq:ode}, which manifests itself via the right-hand-side $
\f(\x, \aalpha):\R^d\times \R^\ell\rightarrow\R^d$, is unknown. Our goal is to create an accurate numerical model for the governing equation using data of the state variable $\x$. Although similar to the setting of the recent work of discovering unknown dynamical systems \cite{QinWuXiu2019}, our setting here represents a non-trivial extension. That is, not only the form of the right-hand-side $\f$ is unknown, its associated parameters $\aalpha$ are also unknown.

\subsection{Data}

In order to learn the governing equation, we assume trajectory data of the state variables $\x$ are available.
More specifically, let $N_T$ be the total number of trajectories. For each $i$-th trajectory, we have data in the following form,
$$
\mathbf{X}^{(i)} = \left\{\x\left(t_k^{(i)}; \aalpha^{(i)}, \x_0^{(i)}\right)\right\}, \qquad i=1,\dots, N_T, \quad k=0,\dots, K^{(i)},
$$
where $\{t_k^{(i)}\}$ are the time instances where the data are made available,  and the parameter $\aalpha^{(i)}$ and the initial condition $\x_0^{(i)}$ associated with the $i$-th trajectory are unknown.

A distinct feature of the learning method in this paper is to approximate the underlying flow map of the unknown governing equation \eqref{eq:ode}.
The method requires the use of trajectory data from two different time instances. Consequently, we re-organize the data set into pairs of two
adjacent time instances, for each $i=1,\dots, N_T$,
$$
\left\{\x\left(t_k^{(i)}; \aalpha^{(i)}, \x_0^{(i)}\right), \quad\x\left(t_{k+1}^{(i)}; \aalpha^{(i)}, \x_0^{(i)}\right)\right\}, \qquad k=0,\dots, K^{(i)}-1.
$$
Note that for the autonomous system \eqref{eq:ode} considered in this paper, only the time difference is important in the pairs. The actual time value is not
relevant for it can be arbitrarily shifted. Therefore, we write the data pairs as, for each $i=1,\dots, N_T$,
$$
\left\{\x\left(0; \aalpha^{(i)}, \x_0^{(i)}\right), \quad\x\left(\Delta_{k}^{(i)}; \aalpha^{(i)}, \x_0^{(i)}\right)\right\}, \qquad k=0,\dots, K^{(i)}-1,
$$
where $\Delta_k^{(i)} = t_{k+1}^{(i)} - t_{k}^{(i)}$. Finally, to account for possible noises in the data and by using a single index to simplify notation, we
write the data set as
\begin{equation} \label{data_set}
  S = \left\{\z_j^{(1)}, \z_j^{(2)}\right\}, \qquad j=1,\dots, J,
  \end{equation}
  where $J= K^{(1)}+\cdots + K^{(N_T)}$ is the total number of data pairs and
\begin{equation}
\label{eq:data_pair}
\z^{(1)}_j=\left(\x\left(0; \aalpha^{(j)}, \x_0^{(j)}\right)+\epsilon_j^{(1)}, \aalpha^{(j)}, \Delta_j\right) \qquad \z^{(2)}_j=\x\left(\Delta_j; \aalpha^{(j)}, \x_0^{(j)}\right)+\epsilon_j^{(2)},
\end{equation}
where $\epsilon_j^{(1)}$ and $\epsilon_j^{(2)}$ are noises/errors in the state variable data. That is, each $j$-th pair, $j=1,\dots, J$, consists of data of the state variables $\x$
separated by a time difference $\Delta_j$. The pair resides on a certain trajectory associated with an unknown parameter value $\aalpha^{(j)}$ and is
originated from an unknown initial condition $\x_0^{(j)}$.
We also define
\begin{equation} \label{I_delta}
  I_\Delta = [\min_j \Delta_j, \max_j \Delta_j]
\end{equation}
to be the range of the time lags in the dataset \eqref{data_set}.

\subsection{Deep neural networks}\label{sec:FNN}

In this paper we adopt deep neural network (DNN) as the primary modeling method for recovering unknown governing equation. In particular, we
employ feed forward neural networks (FNN) as the core building block. A standard FNN defines a nonlinear map in the following sense.
Let $\N:\R^m\rightarrow\R^n$ be the operator associated with a FNN with $M\geq 1$ hidden layers. Its map can be written as
\begin{equation}
\label{eq:fnn}
\y^{out}=\N(\y^{in};\Theta)=\W_{M+1}\circ(\sigma_M\circ \W_{M})\circ\cdots\circ (\sigma_1\circ \W_1) (\y^{in}),
\end{equation}
where $\W_j$ is weight matrix between the $j$-th layer and the $(j+1)$-th layer, $\sigma_j:\R\rightarrow \R$ is the activation function, and $\circ$ stands for operator composition. Following the standard notation, we have augmented biases into the weight matrices, and the activation function is applied in component-wise manner.  To simplify the notation we use $\Theta$ to denote all the model parameters $\Theta=\{\W_j\}_{j=1}^{M+1}$ in the FNN.

When the input and output dimensions are identical, i.e., $m=n$, residue network (ResNet) can be readily defined as
\begin{equation} \label{ResNet}
  \y^{out}=\left[\mathbf{I}_m + \N(\cdot;\Theta)\right](\y^{in}),
\end{equation}
where $\mathbf{I}_m$ is the identity matrix of size $m\times m$. In this form, the neural network in fact models the difference between the input and output (thus the term ``residue''). Although mathematically equivalent to the original standard DNN, ResNet has been shown to be exceptionally useful in practice after
its introduction in \cite{he2016deep}. We will adopt the ResNet idea and modify it to our modeling work in the following section.

\section{Main Method}\label{sec:method}

In this section we present the main method for recovering
\eqref{eq:ode}. We first present our neural network approximation for
the unknown equation in the context of flow map modeling. We then
discuss its application for uncertainty quantification. Finally, we
present analysis on the error estimation of our method.

\subsection{Parameterized Flow Map}

For given initial condition $\x_0\in I_\x$ and parameter $\aalpha\in I_{\aalpha}$, the (unknown) dynamical system \eqref{eq:ode} defines a flow map
\begin{equation}
\label{eq:flow_map}
\x(s; \aalpha, \x_0)=\PPhi_{s-s_0}(\x(s_0; \aalpha, \x_0), \aalpha), 
\end{equation}
which maps the solution $\x$ at time $s_0$ to the one at another different time
$s$. Once again, only the time difference $s-s_0$ is relevant for
autonomous systems considered in this paper.
Consequently, for given time difference $\delta\in I_{\Delta}$, $\PPhi_\delta:\mathbb{R}^d\times\mathbb{R}^l\to\mathbb{R}^d$ defines a map
such that, for any parameter $\aalpha\in I_{\aalpha}$ and initial
condition $\x_0\in I_{\x}$,
\begin{equation}
\label{eq:flow_map_Delta_lag}
	\x(\delta; \aalpha, \x_0)=\PPhi_\delta(\x_0, \aalpha).
\end{equation}
For notational clarity, in the following we suppress the dependence on
the initial condition $\x_0$ and use $\x(\delta;\aalpha)$ in place of $\x(\delta;\aalpha, \x_0)$.

Next, let us characterize the structure of the flow map more carefully. If we integrate the ODE \eqref{eq:ode} from $0$ to $\delta$, we have
\begin{equation}
\label{eq:x_delta}
\x(\delta; \aalpha)=\x(0; \aalpha)+\int_{0}^{\delta} \f(\x(s;
\aalpha), \aalpha)\, ds=\x(0; \aalpha)+\int_{0}^{\delta}
\f(\PPhi_s(\x(0; \aalpha)), \aalpha) ds.
\end{equation} 
We then have, for any fixed $\aalpha$ and $\delta$,
\begin{equation} \label{fmap}
  \x(\delta; \aalpha)=\left[\mathbf{I}_d + {\bf{\Psi}}(\cdot, \aalpha, \delta)\right](\x(0; \aalpha)),
\end{equation}
where $\mathbf{I}_d$ is the identity matrix of size $d\times d$, and for any $\z\in I_{\x}$,
\begin{equation} \label{Psi}
{\bf{\Psi}}(\z, \aalpha, \delta) = \int_{0}^{\delta}
\f\left(\PPhi_s(\z, \aalpha), \aalpha\right) ds
\end{equation}
is the effective increment along a trajectory from $\z$ over time lag
$\delta$. Upon comparing \eqref{fmap}  with the ResNet operator
\eqref{ResNet}, we see that ResNet provides a natural representation
for flow map increment for fixed time lag $\delta$ and parameter
$\aalpha$. 
We remark that this is
an exact representation and without any approximation over the time
horizon $\delta$. Therefore, ResNet structure can be considered an
exact time integrator in term of flow map (\cite{QinWuXiu2019}).

\subsection{Neural Network Model Construction}

\begin{figure}
	\centering
	\includegraphics[scale=0.4]{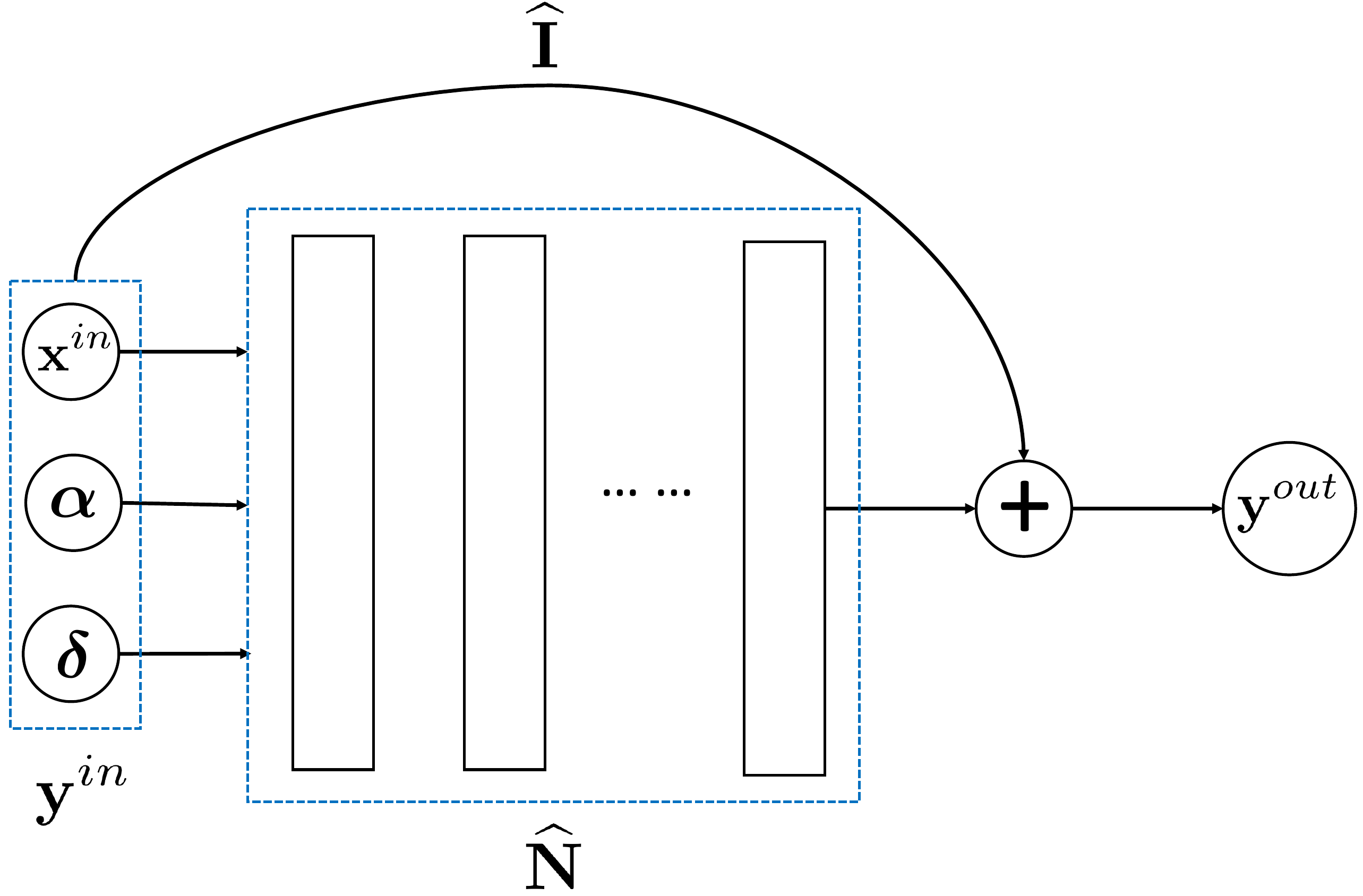}
	\caption{Structure of the neural network.}
	\label{fig:ResNet}
      \end{figure}

Straightforward application of ResNet, as proposed in
\cite{QinWuXiu2019}, is only applicable for \eqref{fmap} for fixed
$\aalpha$ and $\delta$. To model unknown system \eqref{eq:ode} with
unknown parameters, we therefore propose a modified ResNet
structure. The network structure is illustrated in
Fig.~\ref{fig:ResNet}. The input consists of the state
variable $\x^{in}$, the parameters $\aalpha$ and the time lag $\delta$,
i.e., $\y^{in} = [\x^{in}, \aalpha, \delta]^\top$.
Let us
define a $(d+\ell+1)\times (d+\ell+1)$ matrix
\begin{equation}
  \widehat{\mathbf{I}} = \left[
      \begin{array}{cc}
        \mathbf{I}_d & \mathbf{0} \\
        \mathbf{0} & \mathbf{0}
      \end{array}
      \right],
    \end{equation}
    where $\mathbf{I}_d$ is the identity matrix of size $d\times
    d$. Then the network operator is defined as
    \begin{equation}
      \y^{out} =  \left[\widehat{\mathbf{I}} +
        \widehat{\N}\right]\left(\y^{in}\right),
    \end{equation}
    where $\hN:\R^{d+\ell+1}\to\R^d$ is the operator corresponding to the deep FNN. It is obvious that this is equivalent to
    \begin{equation}
      \x^{out}(\x^{in}, \aalpha, \delta; \Theta) = \x^{in} +
      \widehat{\N}(\x^{in}, \aalpha, \delta; \Theta),
    \end{equation}
where $\Theta$ stands for the set of parameters in the neural network structure.
    By using the pair-wise data set \eqref{data_set}, we set
    $$
    \z_j^{(1)}\rightarrow (\x^{in}, \aalpha, \delta), \qquad \z_j^{(2)}\rightarrow \x^{out}.
    $$
    The network training is conducted to find the neural network
    parameter set $\Theta^*$ that minimizes the mean-squared loss, i.e.,
    \begin{equation} \label{eq:loss}
      \Theta^* =\argmin_\Theta\frac{1}{J}\sum_{j=1}^J\left\|
        \x^{out}(\z_j^{(1)}; \Theta)-\z_j^{(2)}\right\|^2,
      \end{equation}
    Upon finding the optimal network parameter $\Theta^*$, we obtain trained network model
    \begin{equation} \label{model}
      \x(\delta; \aalpha) = \x(0; \aalpha) + \widehat{\N}(\x(0;
      \aalpha), \aalpha, \delta; \Theta^*), \qquad \aalpha\in I_{\aalpha}, \quad
      \delta\in I_\Delta,
    \end{equation}
    where $I_\Delta$ is the range
    of time lags in the dataset \eqref{data_set} and is defined in \eqref{I_delta}.
    Upon comparing to \eqref{fmap}, it is obvious that the trained
    neural network $\widehat{\N}$ is an approximation to the parameterized effective
    increment ${\bf \Psi}$ in \eqref{Psi}. That is,
    \begin{equation} \label{N_approx}
      \widehat{\N}(\z, \aalpha, \delta) \approx {\bf{\Psi}}(\z,
      \aalpha, \delta), \qquad (\z,\aalpha,\delta)\in I_{\x}\times
      I_{\aalpha} \times I_\Delta.
    \end{equation}
    The error in this approximation depends on
    the data quality and numerical training algorithm. Despite of this
    approximation error, the trained model \eqref{model} does not
    contain any error associated with time discretization. Therefore,
    the model \eqref{model} can be considered as an exact time integrator
    for the unknown parameterized system \eqref{eq:ode}.


    \subsection{Prediction and Uncertainty Quantification}

    Once the model \eqref{model} is constructed, it serves as an
    approximation to the flow map of the original system
    \eqref{eq:ode}. We can conduct system prediction via iterative use
    of the model. Let $\delta_k\in I_\delta$ be a sequence of time
    steps and $\x_0$ be a given initial condition. Then, for given
    system parameter $\aalpha\in I_{\aalpha}$, we have
    \begin{equation} \label{prediction}
      \left\{
      \begin{split}
        &\wh{\x}(t_0; \aalpha) = \x_0, \\
      &\wh{\x}(t_{k+1}; \aalpha) = \wh{\x}(t_{k}; \aalpha) + \widehat{\N}(\wh{\x}(t_k;
        \aalpha), \aalpha, \delta_k; \Theta^*), \\
        &t_{k+1} =
        t_k + \delta_k, \qquad k=0,1,\dots.
      \end{split}
      \right.
    \end{equation}
These serve as approximation of the true solution $\x(t; \aalpha,
\x_0)$ of the unknown system \eqref{eq:ode} at the time instances
$t\in\{t_k, k=0,1,\dots\}$, with given parameter value $\aalpha$ and initial
condition $\x_0$.

    When the system parameters $\aalpha$ are random, uncertainty
    quantification can be readily conducted by using the system
    prediction model \eqref{prediction}.
Let $\rho_{\aalpha}$ be the probability density the parameters
$\aalpha$. Statistical information of the true solution $\x(t; \aalpha)$  of \eqref{eq:ode} can approximated by applying the required
statistical analysis on the approximate solution \eqref{prediction} generated by the
network model \eqref{model}, as a post-processing step.
For example, the mean and variance of the solution 
can be approximated as
\begin{equation}
\label{eq:mean_var_approx}
\begin{split}
\Eax&\approx \Eahx= \int_{I_{\aalpha}} \wh{\x}(t; \y) \rho_{\aalpha} (\y)\, d\y, \\
\Vax&\approx \Vahx = \int_{I_{\aalpha}}\left[\wh{\x}(t; \y)-\E_{\aalpha}[\wh{\x}(t; \y)]\right]^2\rho_{\aalpha}(\y)\,d\y. 
\end{split}
\end{equation}
The integrals of $\wh{\x}$ can be further approximated sampling based
method, e.g., Monte Carlo or quadratue rule. This requires
the predictive solution $\wh{\x}$ at different sampling locations of
$\aalpha$, which can be produced by the learned network model \eqref{prediction}.

\subsection{Theoretical properties} \label{sec:analysis}

We now establish some theoretical analysis of the proposed
methods.
Our analysis relies on the celebrated universal approximation result for single hidden layer
full connected feedforward neural network.
\begin{theorem} {\cite[Theorem 3.1]{pinkus1999}}
	For any given function $f\in C(\mathbb{R}^n)$ and positive real number $\epsilon>0$, there exists a single-layer neural network $N(\cdot\,; \Theta)$ with parameter $\Theta$ such that 
	\begin{equation*}
		\max_{\x\in D} |f(\x)-N(\x\,;\Theta)| <\epsilon
	\end{equation*}
	for any compact set $D\in \mathbb{R}^n$, if and only if the activation functions are continuous and are not polynomials.
\end{theorem}

Based on this result, we assume that the deep neural network used in
our method can be trained such that the error in the approximation
\eqref{N_approx} can be sufficiently small. That is,
for a given $\mathcal{E}>0$
\begin{equation}
  \label{lem:assumption}
  \left|\widehat{\N}(\z, \aalpha, \delta; \Theta^*) - {\bf{\Psi}}(\z,
      \aalpha, \delta)\right| \leq  \mathcal{E} , \qquad \forall (\z,\aalpha,\delta)\in I_{\x}\times
      I_{\aalpha} \times I_\Delta.
\end{equation}
In the following, for any given $n$,  set $t=\sum_{i=0}^{n-1}\delta_i$ as final time for the prediction, where $\delta_i\in I_{\Delta}$ as defined in \eqref{prediction}.
\begin{lemma}
	\label{lem}
	Suppose the right-hand-side $\f(\x, \aalpha)$ of
        \eqref{eq:ode} is Lipschitz continuous with respect to $\x$
        with a uniform Lipschitz constant $L$, for all $\aalpha\in I_{\aalpha}$. If the trained neural network satisfies \eqref{lem:assumption}, then 
	\begin{equation}
	\label{eq:approx_err}
		\norm{\wh{\x}(t;\, \cdot\,,\,\cdot)-\x(t; \,\cdot\,,\,\cdot)}_{L^\infty(I_{\aalpha}\times I_{\x})}\leq \frac{e^{nL\Delta}-1}{e^{L\Delta}-1} \, \mathcal{E},
	\end{equation}
	where $\x(t; \aalpha, \x_0)$ and $\wh{\x}(t; \aalpha, \x_0)$
        are as defined in \eqref{eq:flow_map} and \eqref{prediction}, respectively.
\end{lemma}
\begin{proof}
First, for $(\v, \aalpha, \delta)\in I_\x\times I_{\aalpha}\times I_\delta$, by combining \eqref{eq:flow_map_Delta_lag} and \eqref{fmap}, we have
\begin{equation*}
	\PPhi_\delta (\v, \aalpha)=\v+\PPsi(\v, \aalpha, \delta).
\end{equation*}
Furthermore, \eqref{model} defines an approximated flow map 
\begin{equation*}
	\wh{\PPhi}_\delta(\v, \aalpha;\Theta^*)=\v+\hN(\v, \aalpha, \delta; \Theta^*), \quad (\v, \aalpha, \delta)\in I_\x\times I_{\aalpha}\times I_\delta.
\end{equation*}
Next, for any $\u, \v\in I_\x$, and fixed $(\aalpha, \delta)\in I_{\aalpha}\times I_\Delta$, let us consider
\begin{align}
\label{lem:eq}
\nonumber
 \abs{\PPhi_\delta(\u, \aalpha)-\wh{\PPhi}_\delta(\v, \aalpha; \Theta^*)}
 &\leq \abs{\PPhi_\delta(\u, \aalpha)-\PPhi_\delta(\v, \aalpha)}+\abs{\PPhi_\delta(\v, \aalpha)-\wh{\PPhi}_\delta(\v, \aalpha; \Theta^*)}\\\nonumber
 &=\abs{\PPhi_\delta(\u, \aalpha)-\PPhi_\delta(\v, \aalpha)}+\abs{\PPsi(\v, \aalpha, \delta)-\hN(\v, \aalpha, \delta; \Theta^*)}\\
 &\leq e^{L\Delta}|\u-\v|+\mathcal{E}
\end{align}
where in the last step we have used \eqref{lem:assumption} and the classical result on the continuity of dynamical system with respect to the initial data; see \cite[p. 109]{stuart1998dynamical}.

Then for $t=\sum_{i=0}^{n-1} \delta_i$ with $\delta_i\in I_\Delta$, by \eqref{prediction}, we have
\begin{equation*}
\wh{\x}(t; \z, \aalpha)=\wh{\PPhi}_{\delta_{n-1}}\circ\wh{\PPhi}_{\delta_{n-2}}\circ \cdots \circ\wh{\PPhi}_{\delta_0}(\z, \aalpha; \Theta^*).
\end{equation*}
Moreover, by the time invariance of autonomous system, we have
\begin{equation*}
\x(t; \z, \aalpha)={\PPhi}_{\delta_{n-1}}\circ {\PPhi}_{\delta_{n-2}}\circ \cdots \circ {\PPhi}_{\delta_0}(\z; \aalpha).
\end{equation*}
Then, by repeatedly employing \eqref{lem:eq}, we have
\begin{align*}
	&|\wh{\x}(t; \z, \aalpha)-\x(t; \z, \aalpha)| \\
	= & \abs{\wh{\PPhi}_{\delta_{n-1}}\circ\wh{\PPhi}_{\delta_{n-2}}\circ \cdots \circ\wh{\PPhi}_{\delta_0}(\z, \aalpha; \Theta^*)-{\PPhi}_{\delta_{n-1}}\circ {\PPhi}_{\delta_{n-2}}\circ \cdots \circ {\PPhi}_{\delta_0}(\z; \aalpha)}\\
	\leq & \mathcal{E}+e^{L\delta_{n-1}}\abs{\wh{\PPhi}_{\delta_{n-2}}\circ \cdots \circ\wh{\PPhi}_{\delta_0}(\z, \aalpha;\Theta^*)- {\PPhi}_{\delta_{n-2}}\circ \cdots \circ {\PPhi}_{\delta_0}(\z;\aalpha)}\\
	\leq & \mathcal{E}+e^{L\delta_{n-1}}\left[\mathcal{E}+e^{L\delta_{n-2}}\abs{\wh{\PPhi}_{\delta_{n-3}}\circ \cdots \circ\wh{\PPhi}_{\delta_0}(\z, \aalpha; \Theta^*)-{\PPhi}_{\delta_{n-3}}\circ \cdots \circ {\PPhi}_{\delta_0}(\z;\aalpha)}\right]\\
	\leq & \ldots\\
	\leq & \mathcal{E}\left(1+e^{L\delta_{n-1}}+e^{L(\delta_{n-1}+\delta_{n-2})}+\ldots+e^{L\sum_{i=1}^{n-1}\delta_i}\right)\\
	\leq & \mathcal{E}\left(1+e^{L\Delta}+e^{2L\Delta}+\ldots+e^{(n-1)L\Delta}\right)\\
	= & \frac{e^{nL\Delta}-1}{e^{L\Delta}-1}\,\mathcal{E}
\end{align*}
for any $(\aalpha, \z)\in I_{\aalpha}\times I_\x$. This implies the result \eqref{eq:approx_err}.
\end{proof}
%

When the trained network is used for UQ as in
\eqref{eq:mean_var_approx},
we obtain the following estimates.
\begin{theorem}
	\label{thm}
	Under the same assumptions of Lemma \ref{lem} and assume for a
        fixed initial condition $\x_0$, the solution is bounded
        $\norm{\x(t;
          \,\cdot)}_{L^\infty(I_{\aalpha})}=C_t<\infty$. Then,
	\begin{align}
	\label{eq:estimate_mean}
			\abs{\Eax-\Eahx}&\leq \,C(n, L, \Delta)\mathcal{E},\\
	\label{eq:estimate_var}
			\abs{\Vax-\Vahx}&\leq 2 \,C(n, L, \Delta)^2\,\mathcal{E}^2+4\,C(n, L, \Delta)\, C_t\, \mathcal{E},
	\end{align}
	where $C(n, L, \Delta)=\frac{e^{nL\Delta}-1}{e^{L\Delta}-1}$.
\end{theorem}

\begin{proof}
	For the mean approximation \eqref{eq:estimate_mean}, by Lemma \ref{lem}, we have
	\begin{equation*}
		\abs{\Eax-\Eahx}\leq \Ea[\abs{\xta-\hxta}]\leq \frac{e^{nL\Delta}-1}{e^{L\Delta}-1}\,\mathcal{E}.
	\end{equation*}
	For the variance approximation, we have
	\begin{align*}
		&\abs{\Vax-\Vahx}\\
		=& \abs{\Ea[\xta^2]-\Eax^2-\Ea[\hxta^2]+\Eahx^2}\\	
		\leq & \abs{\Ea[\xta^2]-\Ea[\hxta^2]}+\abs{\Eax^2-\Eahx^2}\\
		=& \Ea[\abs{\xta-\hxta}\abs{\xta+\hxta}]\\
		&+\Ea[\abs{\xta-\hxta}]\,\Ea[\abs{\xta+\hxta}]\\
		\leq& 2\frac{(e^{nL\Delta}-1)\mathcal{E}}{e^{L\Delta}-1}\, \Ea[\abs{\xta+\hxta}]\\
		\leq& 2\frac{(e^{nL\Delta}-1)\mathcal{E}}{e^{L\Delta}-1}\, (\Ea[\abs{\hxta-\xta}]+\Ea[\abs{2\xta}])\\
		\leq& 2\frac{(e^{nL\Delta}-1)^2}{(e^{L\Delta}-1)^2}\mathcal{E}^2+4\frac{(e^{nL\Delta}-1)C_t}{e^{L\Delta}-1} \mathcal{E}
	\end{align*}
	This gives the result \eqref{eq:estimate_var}.	
\end{proof}

In practice, the true parameter range $\Ia$ is sometimes unknown, due to our lack of knowledge of the underlying physical system and is usually replaced by an estimated range $\tIa$, in which the ODE \eqref{eq:ode} is assumed to be well-posed. Suppose the distribution of parameters on the estimated range $\tIa$ is $\tr$. Then, as in \cite{jakeman2010, chen2013}, the mean and variance can be approximated by
\begin{equation}
\label{eq:mean_var_approx_2}
\begin{split}
\Eax&\approx \widetilde{\E}_{\aalpha}[\hxta]=\int_{\tIa} \wh{\x}(t; \y) \tr (\y)\, d\y, \\
\Vax&\approx \widetilde{\Var}_{\aalpha}[\hxta]=\int_{\tIa}\left[\wh{\x}(t; \y)-\widetilde{E}_{\aalpha}[\wh{\x}(t; \y)]\right]^2\tr(\y)\,d\y.
\end{split}
\end{equation} 
Let us define
$$\Ia^o=\Ia\cap \tIa, \quad \Ia^-=\Ia/\Ia^o,\quad \tIa^-=\tIa/\Ia^o.$$
Suppose 
\begin{equation}
\label{cor:assumption1}
\int_{\Ia^o} \abs{\r(\y)-\tr(\y)}\,d\y \leq \gamma
\end{equation}
and the difference between $\Ia$ and $\tIa$ is small in the sense that 
\begin{equation}
\label{cor:assumption2}
	\int_{\tIa^-} \tr(\y)\,d\y+\int_{\Ia^-} \r(\y)\,d\y \leq \eta
\end{equation}
for some small positive number $\gamma$ and $\eta$, then we have the following results concerning the error in the approximation \eqref{eq:mean_var_approx_2}.
\begin{theorem}
	Under the assumptions in Lemma \ref{lem} and assume the distribution $\tr$ and the estimated parameter range $\tIa$ satisfy \eqref{cor:assumption1}and \eqref{cor:assumption2}. In addition, assume that $\norm{{\x}(t;\,\cdot\,)}_{L^\infty(\tIa\cup\Ia)}\leq \widetilde{C}_t$ 
	, then we have
	\begin{align}
	\label{eq:estimate_mean_2}
	\abs{\Eax-\tEahx}&\leq \,\widetilde{C}_t(\eta+\gamma)+C(1+\eta)\mathcal{E},\\
	\label{eq:estimate_var_2}\nonumber
	\abs{\Vax-\tVahx}&\leq (3\tC_t^2+C\tC_t\mathcal{E})(\eta+\gamma)\\
	    &+(4\tC_t+2C\mathcal{E})(1+\eta)C\mathcal{E},
	\end{align}
	where $C=\frac{e^{nL\Delta}-1}{e^{L\Delta}-1}$.
\end{theorem}
\begin{proof}
First, let us show the error estimate for the approximation of the expectation in \eqref{eq:estimate_mean_2}. To this end, we have the following sequence of estimates.
\begin{align*}
	 &\abs{\Eax-\tEahx}\\
	=&\abs{\int_{\Ia} \xty\r(\y)\,d\y-\int_{\tIa} \hxty\tr(\y)\,d\y}\\
	=& \left\lvert\int_{\Ia^-}\xty\r(\y)\,d\y+\int_{\Ia^o} \xty\r(\y)\,d\y-\int_{\tIa^-} \hxty\tr(\y)\,d\y- \int_{\Ia^o}\hxty\tr(\y)\,d\y \right\rvert\\
	\leq & \int_{\Ia^-}\abs{\xty}\r(\y)\,d\y+\int_{\tIa^-}\abs{\hxty}\tr(\y)\,d\y + \\
	&\abs{\int_{\Ia^o}\xty(\r(\y)-\tr(\y))\,d\y+\int_{\Ia^o}(\xty-\hxty)\tr(\y)\,d\y}\\
	\leq&\norm{\x(t;\cdot)}_{L^\infty} \int_{\Ia^-}\r(\y)\,d\y+\norm{\wh{\x}(t; \cdot)}_{L^\infty}\int_{\tIa^-}\tr(\y)\,d\y
	+\int_{\Ia^o}\abs{\xty}\abs{\r(\y)-\tr(\y)}\,d\y\\
	&+\int_{\Ia^o}\abs{\xty-\hxty}\tr(\y)\,d\y\\
	\leq &\norm{\x(t;\cdot)}_{L^\infty} \int_{\Ia^-}\r(\y)\,d\y
	+\left(\norm{{\x}(t; \cdot)}_{L^\infty}+\norm{{\x}(t; \cdot)-\wh{\x}(t;\cdot)}_{L^\infty}\right)\int_{\tIa^-}\tr(\y)\,d\y\\
	&+\norm{\x(t;\cdot)}_{L^\infty}\int_{\Ia^o}\abs{\r(\y)-\tr(\y)}\,d\y
	+\norm{\x(t;\cdot)-\wh{\x}(t; \cdot)}_{L^\infty}	\int_{\Ia^o}\tr(\y)\,d\y\\
	\leq& \left[\widetilde{C}_t+C\mathcal{E}\right]\eta +\widetilde{C}_t \gamma + C \mathcal{E}\\
	=&\widetilde{C}_t(\eta+\gamma)+C(1+\eta)\mathcal{E}
\end{align*}
which establishes the result \eqref{eq:estimate_mean_2}. In the last step, we have used \eqref{eq:approx_err}, \eqref{cor:assumption1}, and \eqref{cor:assumption2}. 

Next, we consider the approximation of the variance. By the definition of the variance approximation in \eqref{eq:estimate_var_2}, we have
\begin{align*}
	&\abs{\Vax-\tVahx}\\
	=&\abs{\Ea[\xta^2]-\Ea[\xta]^2-\tEa[\xta^2]+\tEa[\xta]^2}\\
	\leq&\abs{\Ea[\xta^2]-\tEa[\xta^2]}+\abs{\Ea[\xta]^2-\tEa[\xta]^2}
\end{align*}
For the first term, since $\norm{\x(t;\cdot)^2}_{L^\infty( \Ia\cup\tIa)}\leq \widetilde{C}^2_t$ and 
\begin{equation*}
\begin{split}
\norm{\x(t;\cdot)^2-\hx(t;\cdot)^2}_{L^\infty(I_{\x}\times \tIa)} 
&\leq \norm{\x(t; \cdot)-\hx(t; \cdot)}_{L^\infty(\tIa)}\norm{\x(t; \cdot)+\hx(t;\cdot)}_{L^\infty(\tIa)} \\
&\leq C \mathcal{E} \left(C\mathcal{E}+2\widetilde{C}_t \right),
\end{split}
\end{equation*}
by the proof for \eqref{eq:estimate_mean_2} we have
\begin{equation}
\label{eq:thm2_est1}
\abs{\Ea[\xta^2]-\tEa[\xta^2]}\leq \tC_t^2(\eta+\gamma)+(1+\eta)(2\tC_tC\mathcal{E}+C^2\mathcal{E}^2).
\end{equation}
For the second term, since 
\begin{equation*}
	\abs{\Ea[\xta]}\leq \tC_t, \quad \abs{\tEa[\hxta]} \leq (\tC_t+C\mathcal{E}),
\end{equation*}
we have
\begin{align}
\nonumber
&\abs{\Ea[\xta]^2-\tEa[\xta]^2}\\\nonumber
\leq& \abs{\Ea[\xta]-\tEa[\xta]}\abs{\Ea[\xta]+\tEa[\xta]}\\
\label{eq:thm2_est2}
\leq&\left[ \tC_t(\eta+\gamma)+C(1+\eta)\mathcal{E}\right](2\tC_t+C\mathcal{E}).
\end{align}
By combining \eqref{eq:thm2_est1} and \eqref{eq:thm2_est2}, we have
\begin{align*}
\abs{\Vax-\tVahx}\leq & (3\tC_t^2+C\tC_t\mathcal{E})(\eta+\gamma)\\
&+(4\tC_t+2C\mathcal{E})(1+\eta)C\mathcal{E}.
\end{align*}
This complete the proof for \eqref{eq:estimate_var_2}.
\end{proof}
\section{Numerical Examples} \label{sec:examples}
In this section we present numerical examples to verify the performance and properties
of the proposed methods.
For benchmarking purpose, we utilize examples with known governing
parameterized equations. We use the true equations only to generate
synthetic data. Once the neural network models are constructed using
the data, we conduct predictions using the trained model and compare
them against the high resolution numerical solutions of the true
equations. In all the examples, time lags, initial
states and parameters are sampled uniformly from $I_\Delta\times
I_{\x}\times I_{\aalpha}$. While we set $I_\Delta = [0, \Delta]$ with
$\Delta = 0.1$ for all examples, the state variable domain $I_{\x}$
and the parameter domain $I_{\aalpha}$ vary and are specified in each example.

%

We use $(M, n)$ to denote the structure of $\hN$, where $M$ denotes
the number of hidden layers and $n$ denotes the number of nodes in
each layer. The activation function is chosen as
$\sigma=\tanh(x)$ (no noticeable difference were observed when using
other activation functions such as ReLU). We
generate data training set by randomly sampling $20 m$ data pairs,
where $m$ is the number of parameters in the network model. The DNNs are trained by minimizing the mean square
loss function in \eqref{eq:loss} 
by using the Adam algorithm \cite{kingma2014adam} with the standard
parameters with the open-source Tensorflow library \cite{tensorflow2015}. The training data set is divided into mini-batches of size $30$. And
we typically train the model for
$2,000$ epochs and reshuffle the training data in each epoch.
All the weights are initialized randomly from Gaussian
distributions and all the biases are initialized to be zeros. 


\subsection{Example 1: Linear Scalar ODE}

Let us first consider the following linear ODE with a single random parameter
\begin{equation}
	\label{eq:example_1}
		\frac{dx}{dt}=-\alpha\, x, \quad x(0)=x_0,
\end{equation}
where $\alpha$ is a random coefficient. We take $I_\alpha=[0,1]$ and $I_x=[0,1]$. 
This simple parameterized equation has the following analytical solution
\begin{equation*}
	x(t;\alpha, x_0)=x_0e^{-\alpha t}.
\end{equation*}

\begin{figure}[htb]
	\centering
	\begin{subfigure}[b]{0.48\textwidth}
		\begin{center}
			\includegraphics[width=1.0\linewidth]{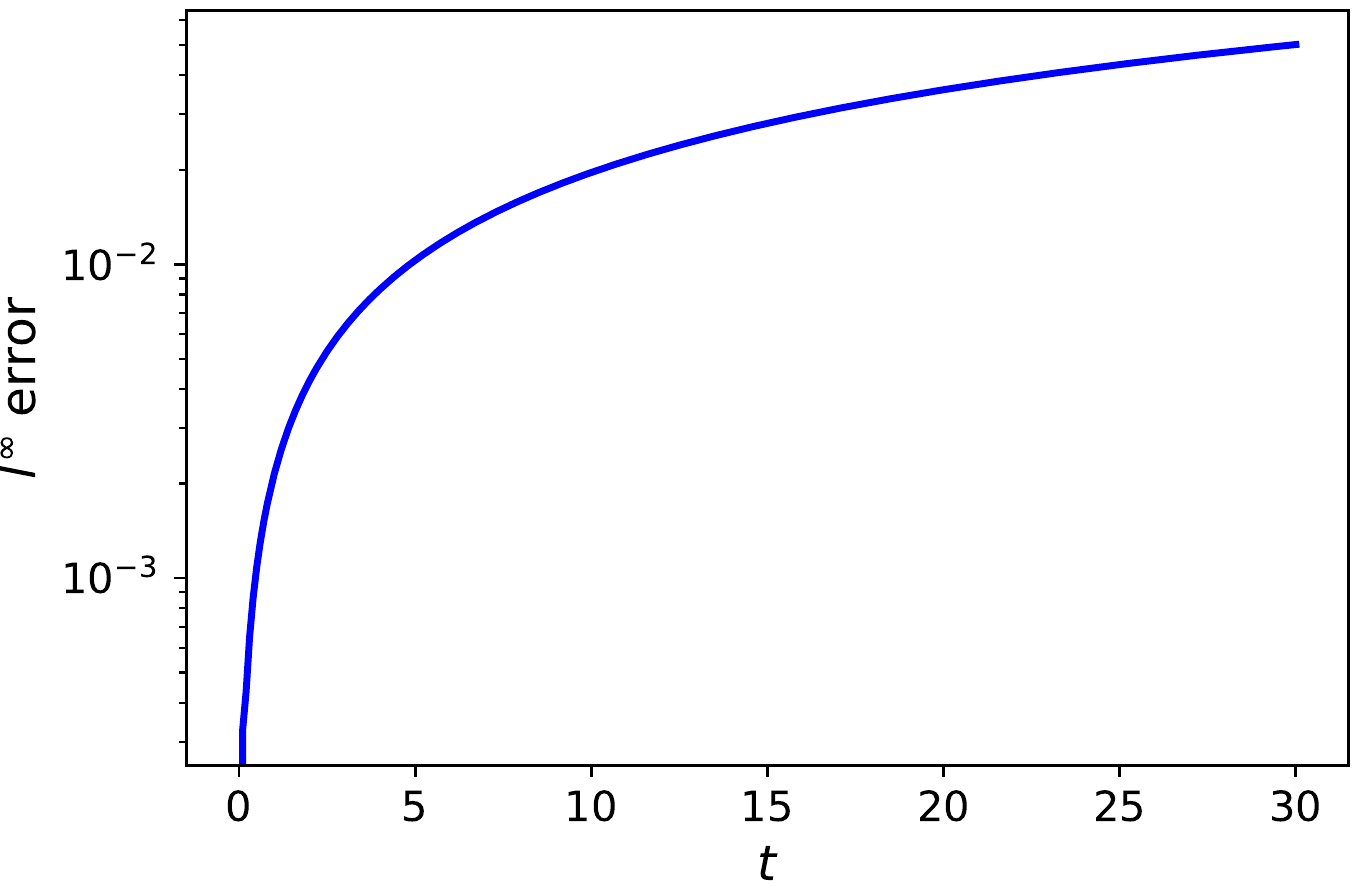}
			\caption{$l^\infty$ error}
		\end{center}
	\end{subfigure}
	\begin{subfigure}[b]{0.48\textwidth}
		\begin{center}
			\includegraphics[width=1.0\linewidth]{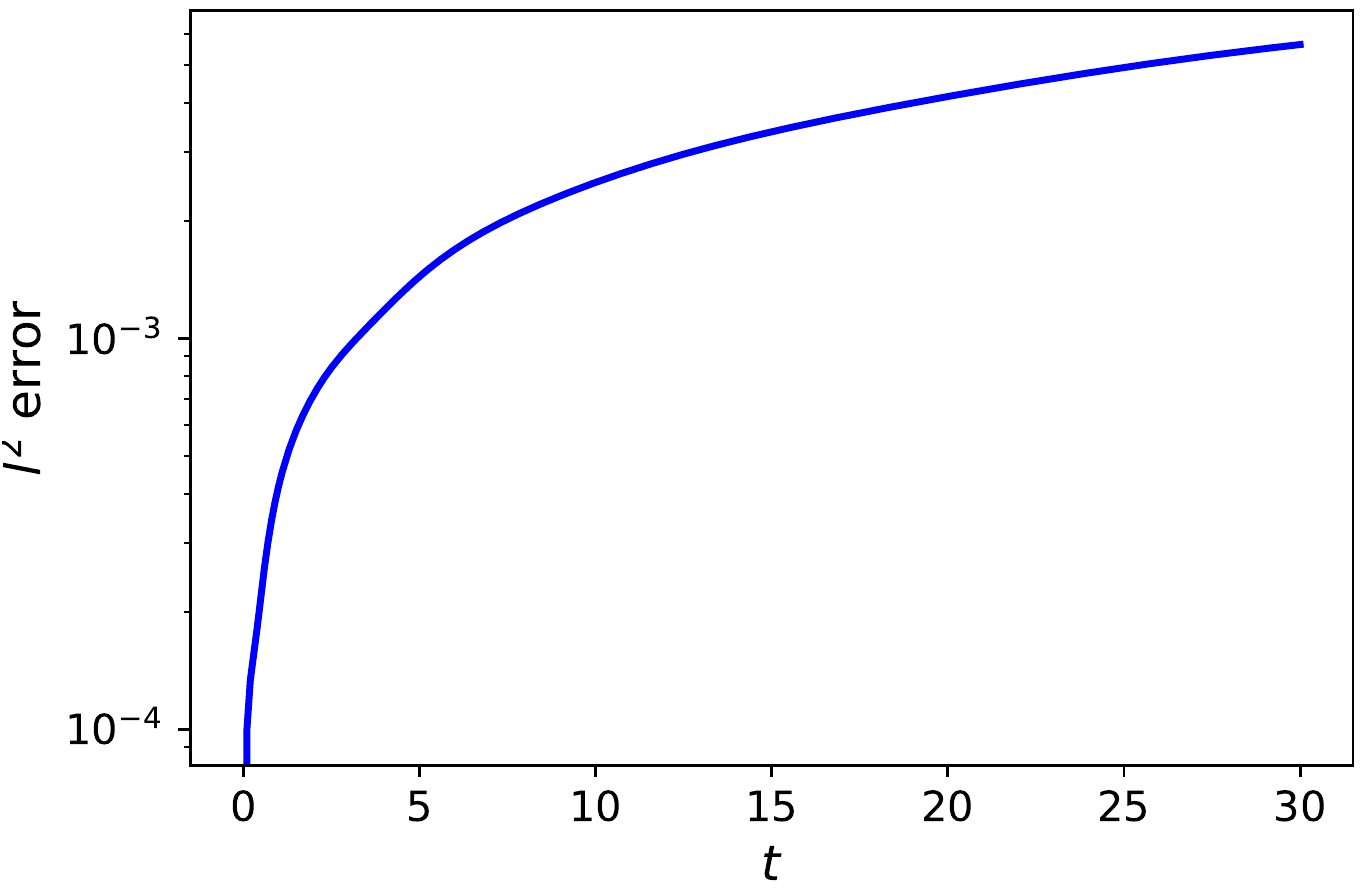}
			\caption{$l^2$ error}
		\end{center}
	\end{subfigure}
	\caption{The $l^\infty$ (left) and $l^2$ (right) error of sample trajectories for Example 1 with $x_0=1$.}
	\label{fig:ex1_error}
\end{figure}

\begin{figure}[htb]
	\centering
	\begin{subfigure}[b]{0.48\textwidth}
		\begin{center}
			\includegraphics[width=1.0\linewidth]{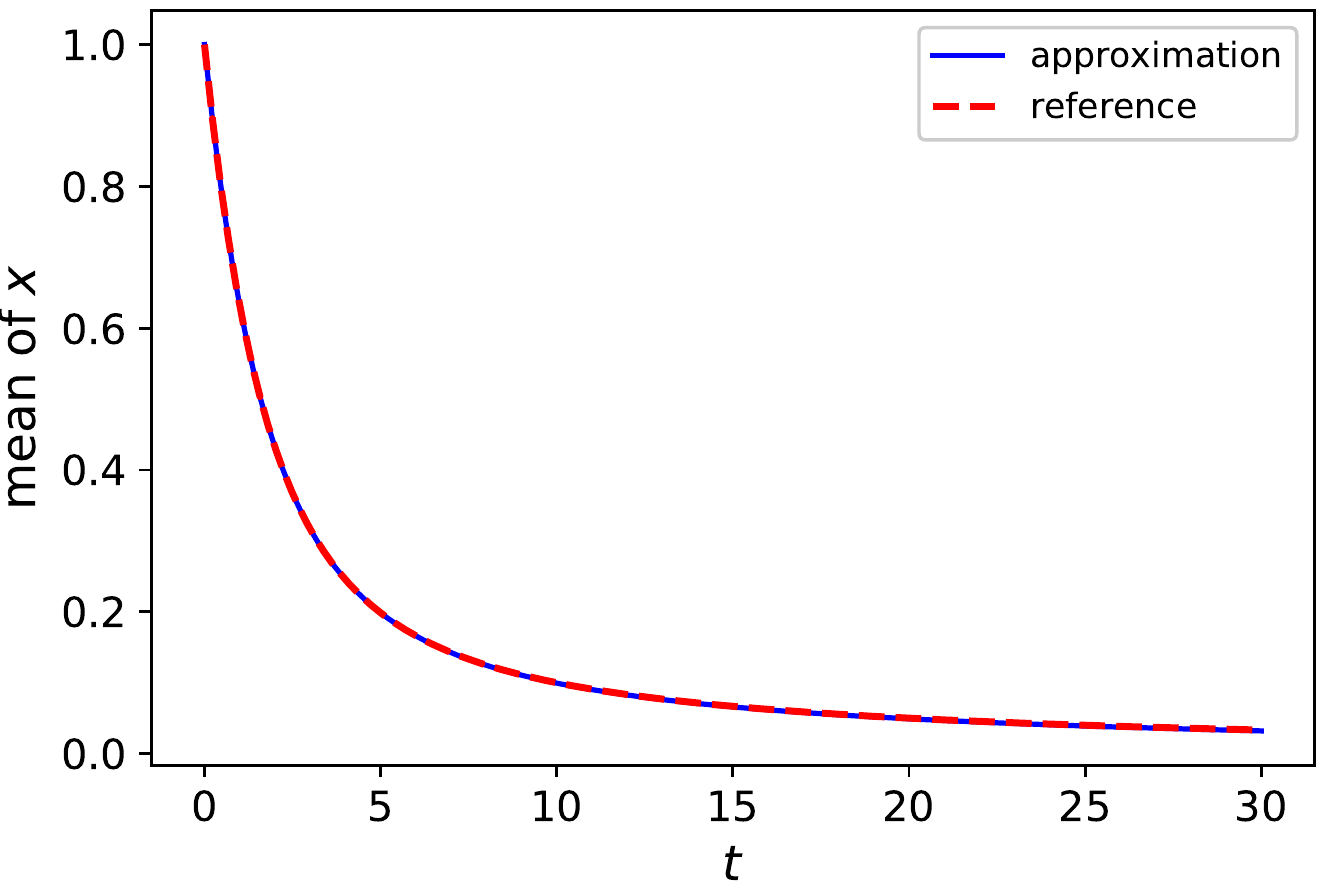}
			\caption{mean of $x(t)$}
		\end{center}
	\end{subfigure}
	\begin{subfigure}[b]{0.48\textwidth}
		\centering
		\includegraphics[width=1.0\linewidth]{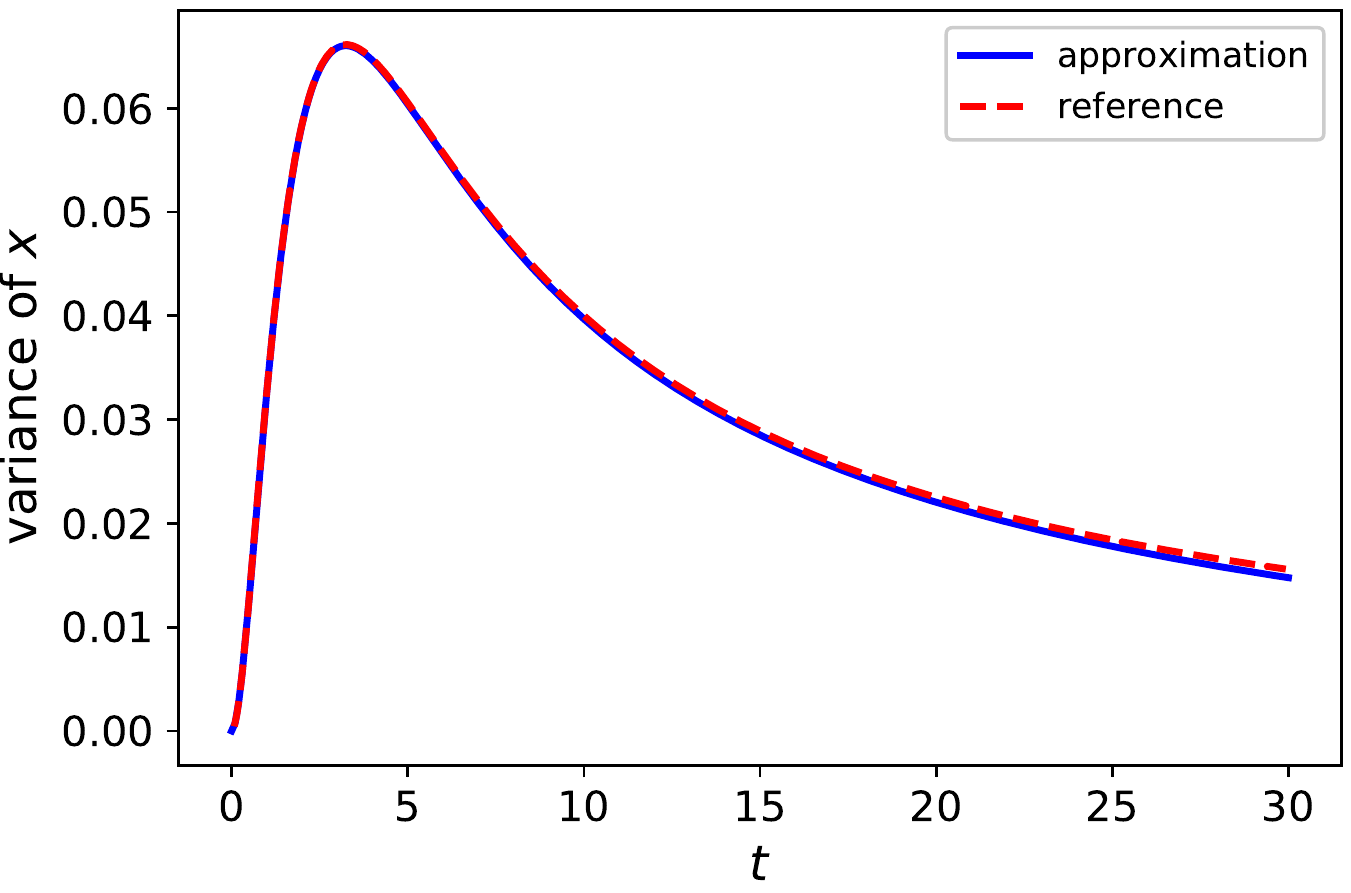}
		\caption{variance of $x(t)$}
	\end{subfigure}	
	\caption{Mean and variance of the solution to Example 1 with $x_0=1$.}
	\label{fig:ex1}
\end{figure}

We approximate the flow map with DNN of structure $(3,40)$. We first
test the $l^\infty$ and $l^2$ error of $100$ sample trajectories. The
trained network with a given sample parameter is composed for $300$
times and makes predictions till $t=30$. The error plots are shown in
\figref{fig:ex1_error}. We observe that the error grows with time,
which is as expected by Lemma \ref{lem}, and stays around $10^{-2}$. 
If $\alpha$ is uniformly distributed in the interval $[0, 1]$, the exact mean of the solution is 
\begin{equation*}
\E_{\alpha}[x(t; \alpha)]=\int_0^1 e^{-\alpha t}\, d\alpha=\frac{1-e^{-t}}{t},
\end{equation*}
and the variance is
\begin{equation}
\Var_{\alpha}[x(t; \alpha)]=\frac{1-e^{-2t}}{2t}-\left( \frac{1-e^{-t}}{t} \right)^2.
\end{equation}
The approximate mean and variance of the DNN approximation of the
governing equations are computed by applying a ten-point
Gauss-Legendre quadrature over the parameter interval $I_\alpha$ to
approximate the integrals in \eqref{eq:mean_var_approx}. In
\figref{fig:ex1}, we show the approximate mean and variance of sample
trajectories. The results are comparable with the results obtained by
the time-dependent gPC in \cite{gerritsma2010}. 
In \figref{fig:ex1_error}, we also present the propagation of errors
in the mean and variance. We observe that the errors in mean solution
grow first and then stays at a level of around $10^{-3}$ after $150$
steps of compositions. The errors in variance continue to grow
exponentially with $n$, as expected by the error estimate in
\eqref{eq:estimate_var}. 
	
\begin{figure}[htb]
	\centering
	\begin{subfigure}[b]{0.48\textwidth}
		\begin{center}
			\includegraphics[width=1.0\linewidth]{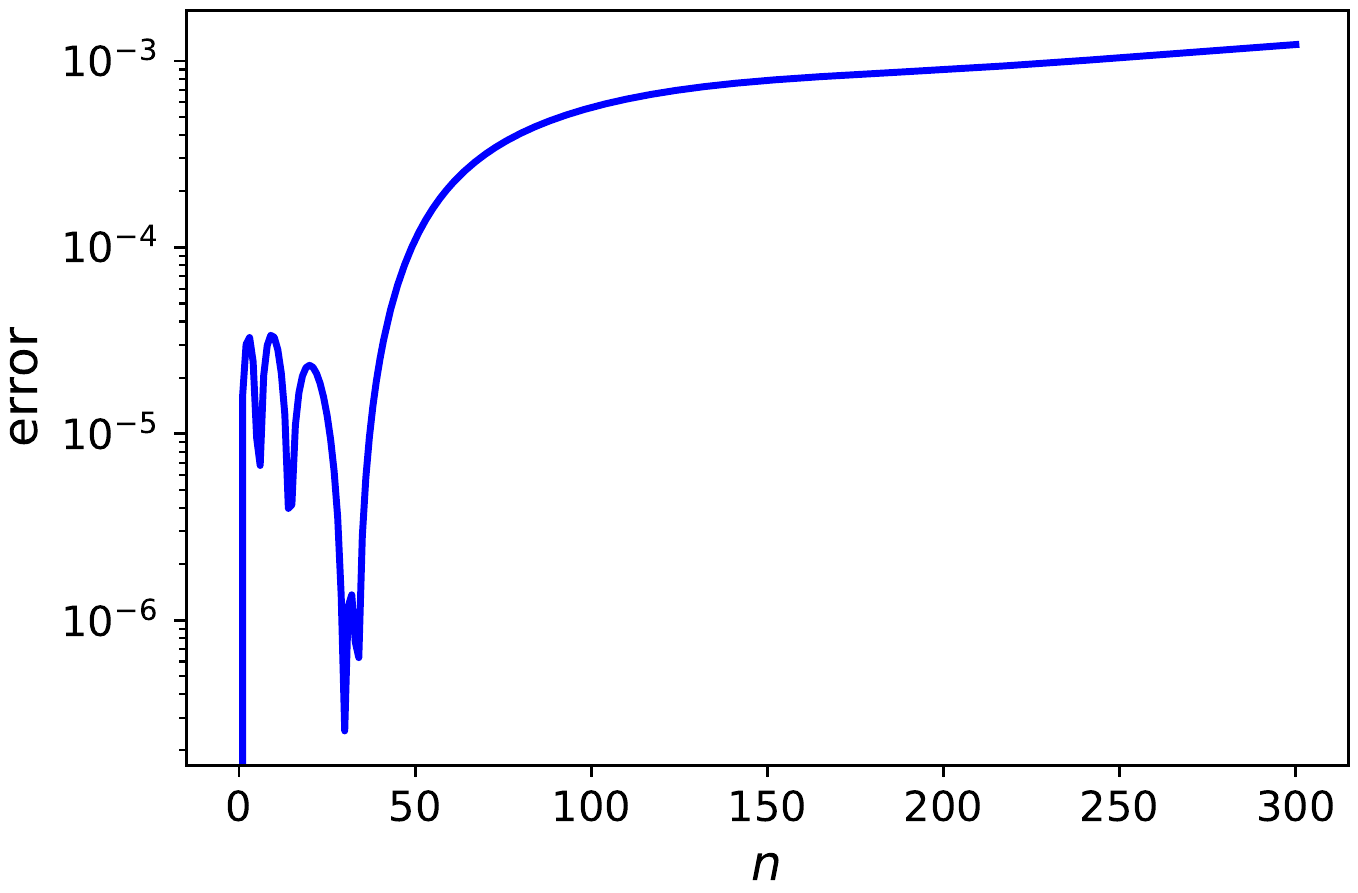}
			\caption{error of the mean}
		\end{center}
	\end{subfigure}
	\begin{subfigure}[b]{0.48\textwidth}
		\centering
		\includegraphics[width=1.0\linewidth]{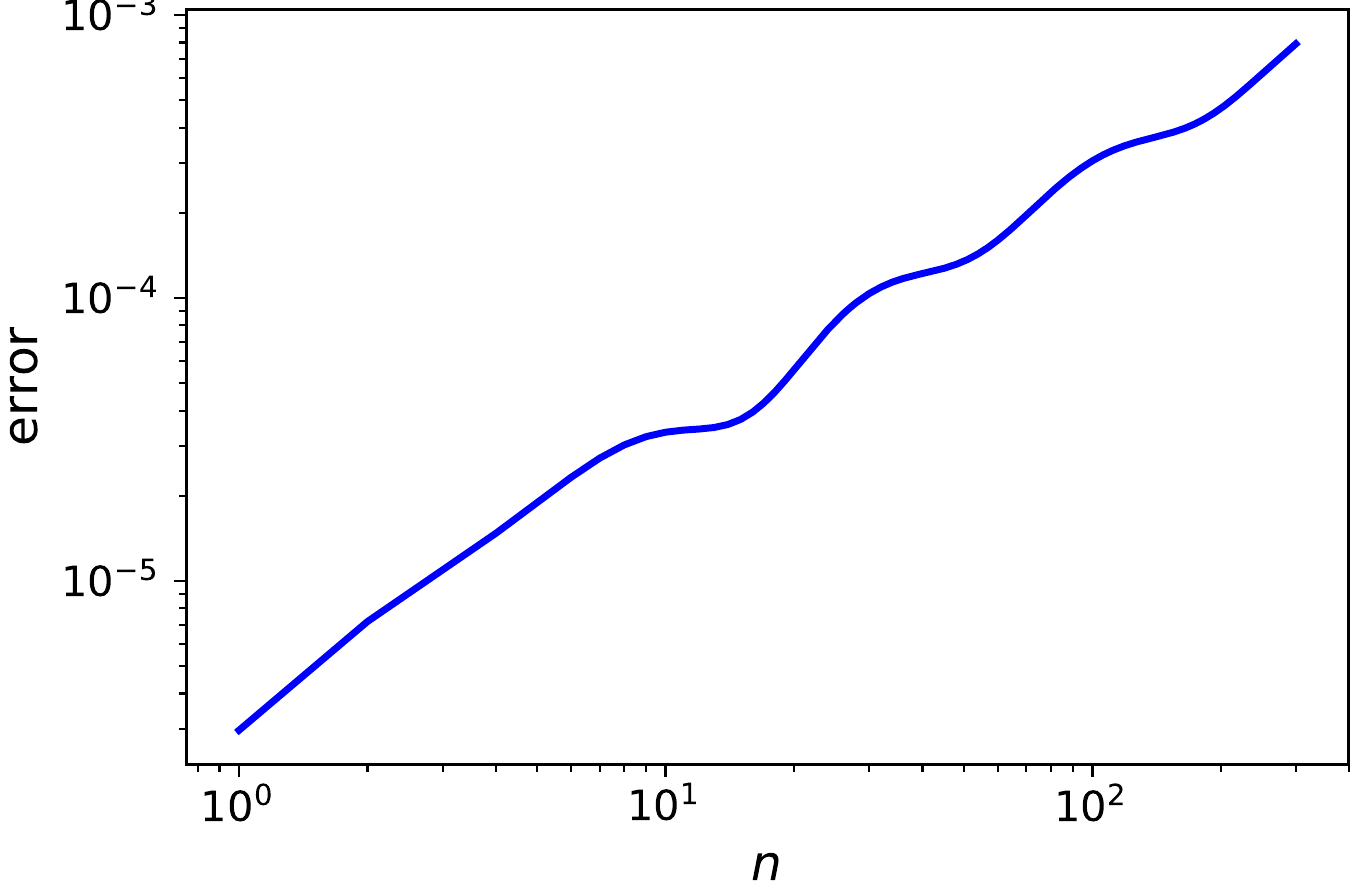}
		\caption{error of the variance}
	\end{subfigure}	
	\caption{The propagation of errors in mean (left)and variance (right) of the solution to Example 1 with $x_0=1$.}
	\label{fig:ex1_error_mean_var}
      \end{figure}

      \subsection{Example 2: Linear ODE System}
      
\begin{figure}[htb]
	\centering
	\begin{subfigure}[b]{0.47\textwidth}
		\begin{center}
			\includegraphics[width=1.0\linewidth]{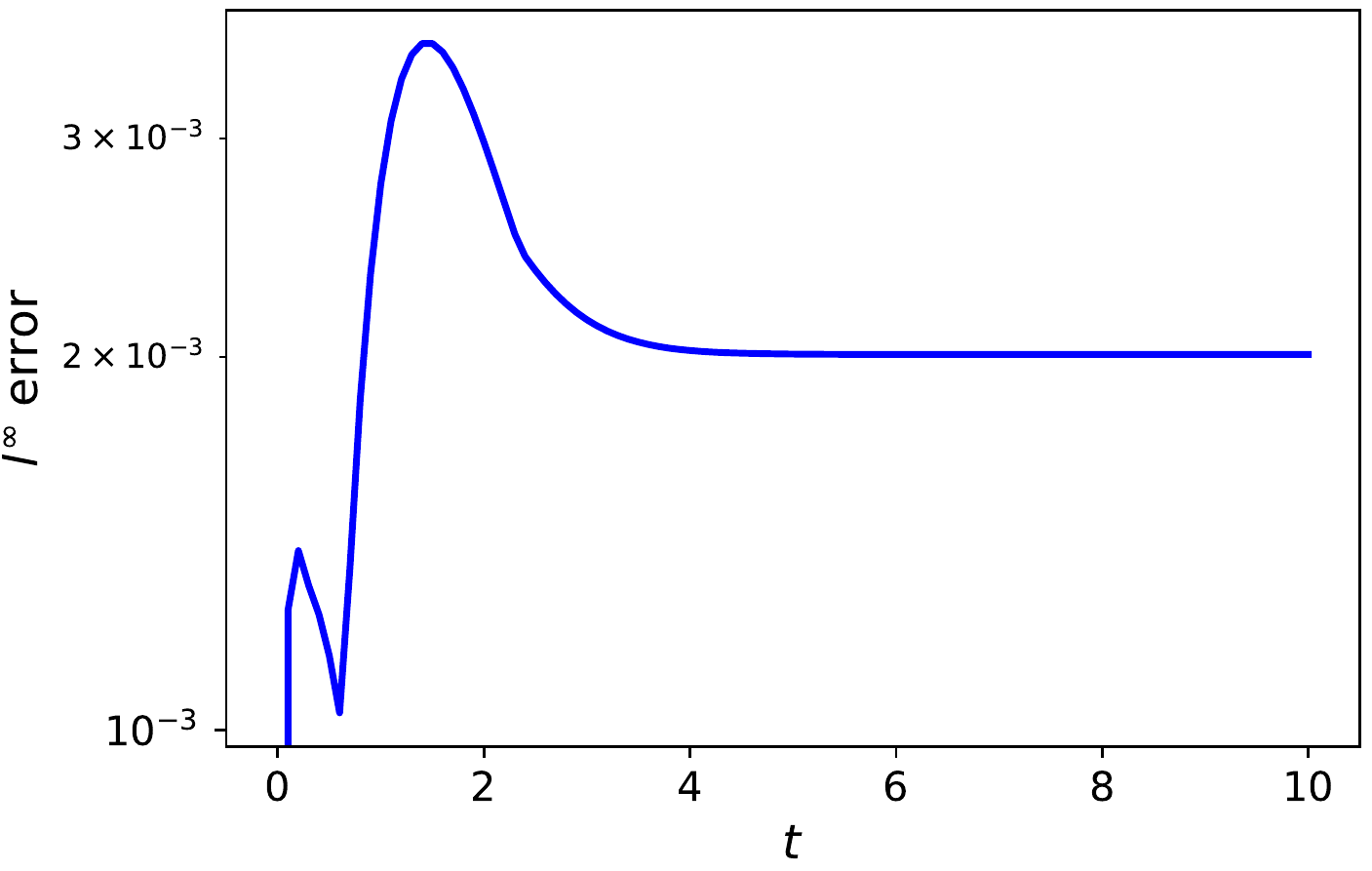}
			\caption{$l^\infty$ error}
		\end{center}
	\end{subfigure}
	\begin{subfigure}[b]{0.47\textwidth}
		\centering
		\includegraphics[width=1.0\linewidth]{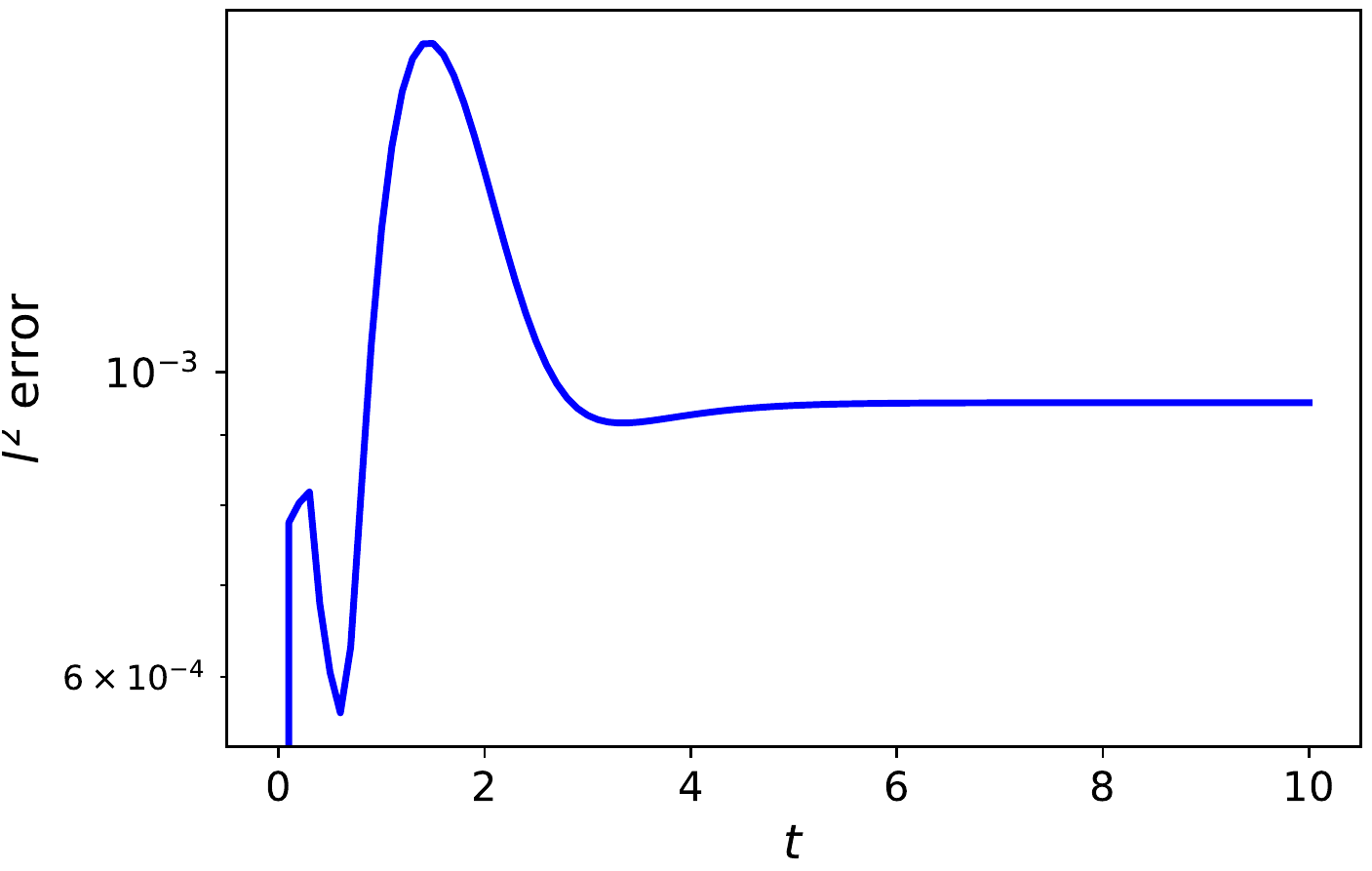}
		\caption{$l^2$ error }
	\end{subfigure}	
	\caption{The $l^\infty$ (left) and $l^2$ (right) error of sample trajectories for Example 2 with $\x_0=(0,1)$.}
	\label{fig:ex2_error}
\end{figure}

We now consider a linear ODE system
\begin{equation}
	\begin{split}
		\frac{dx_1}{dt}&=x_1-\alpha_1 x_2,\\
		\frac{dx_2}{dt}&=\alpha_2 x_1-7x_2,
	\end{split}
\end{equation}
with $\aalpha=(\alpha_1, \alpha_2)\in I_{\aalpha}=(3.8, 4.2)^2$ and
$I_\x=[-1,1]^2$. We use a fully connected block of the structure $(3,
40)$ to construct the DNN model \eqref{model}. After the training of
the DNN is complete, we randomly draw $1,000$ sample parameters
uniformly from $I_{\aalpha}$ and estimate the associated trajectories
for $t\in[0,10]$ and starting from $\x_0=(0,1)$. The $l^\infty$ and
$l^2$ errors for these trajectories are shown in
\figref{fig:ex2_error}. 

\begin{figure}[htb]
	\centering
	\begin{subfigure}[b]{0.48\textwidth}
		\begin{center}
			\includegraphics[width=1.0\linewidth]{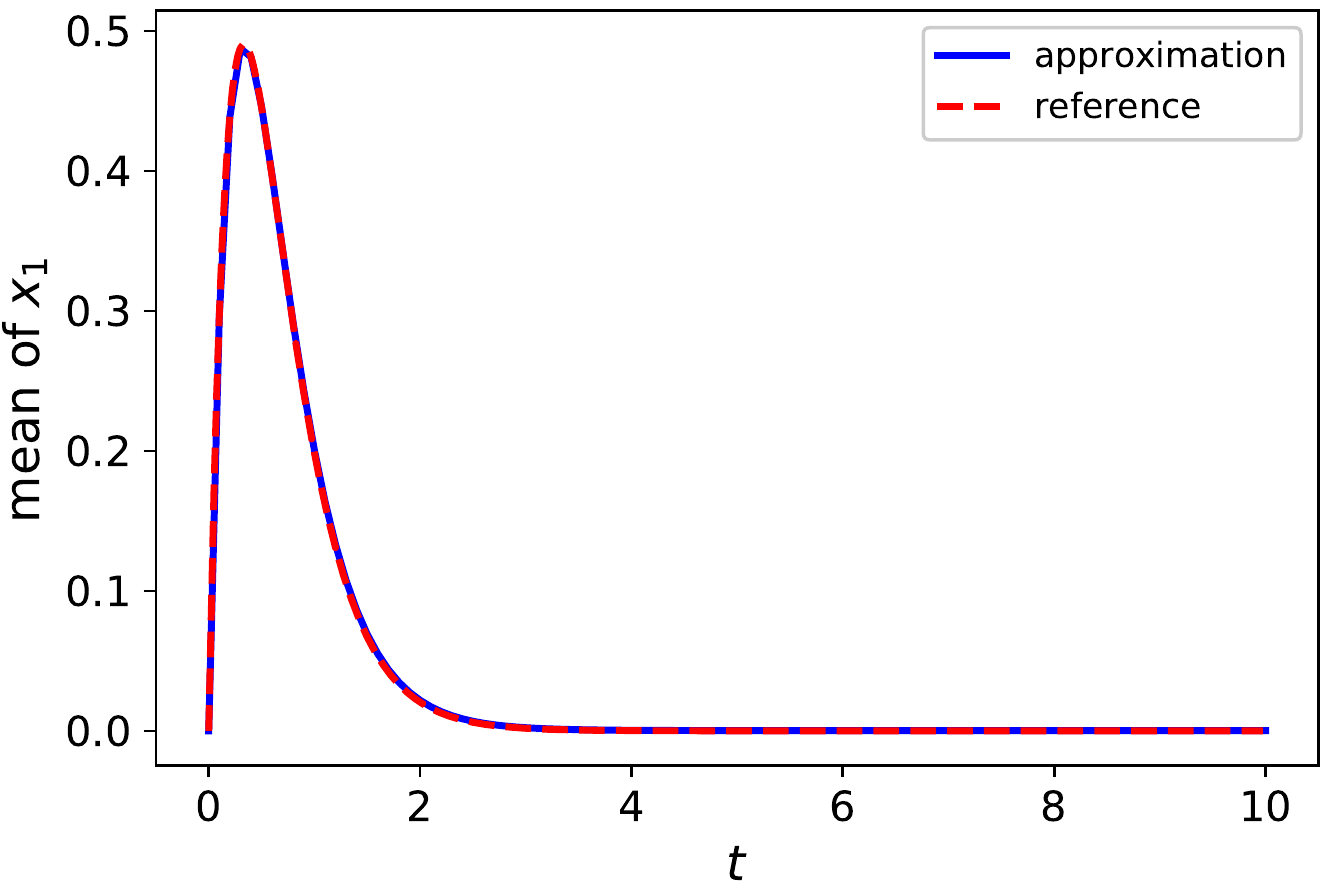}
			\caption{mean of $x_1$}
		\end{center}
	\end{subfigure}
	\begin{subfigure}[b]{0.48\textwidth}
		\centering
		\includegraphics[width=1.0\linewidth]{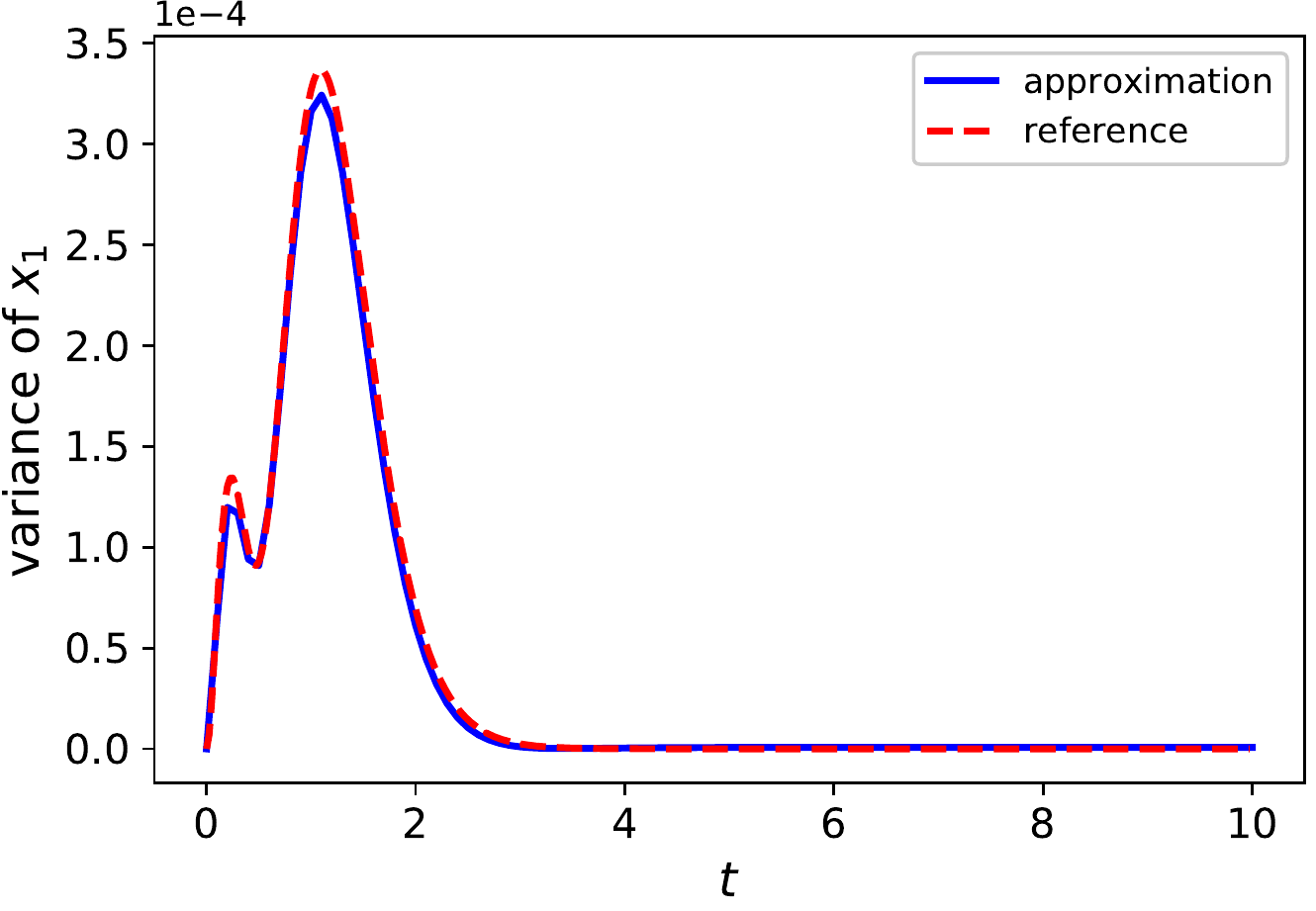}
		\caption{variance of $x_1$}
	\end{subfigure}
	\begin{subfigure}[b]{0.48\textwidth}
		\begin{center}
			\includegraphics[width=1.0\linewidth]{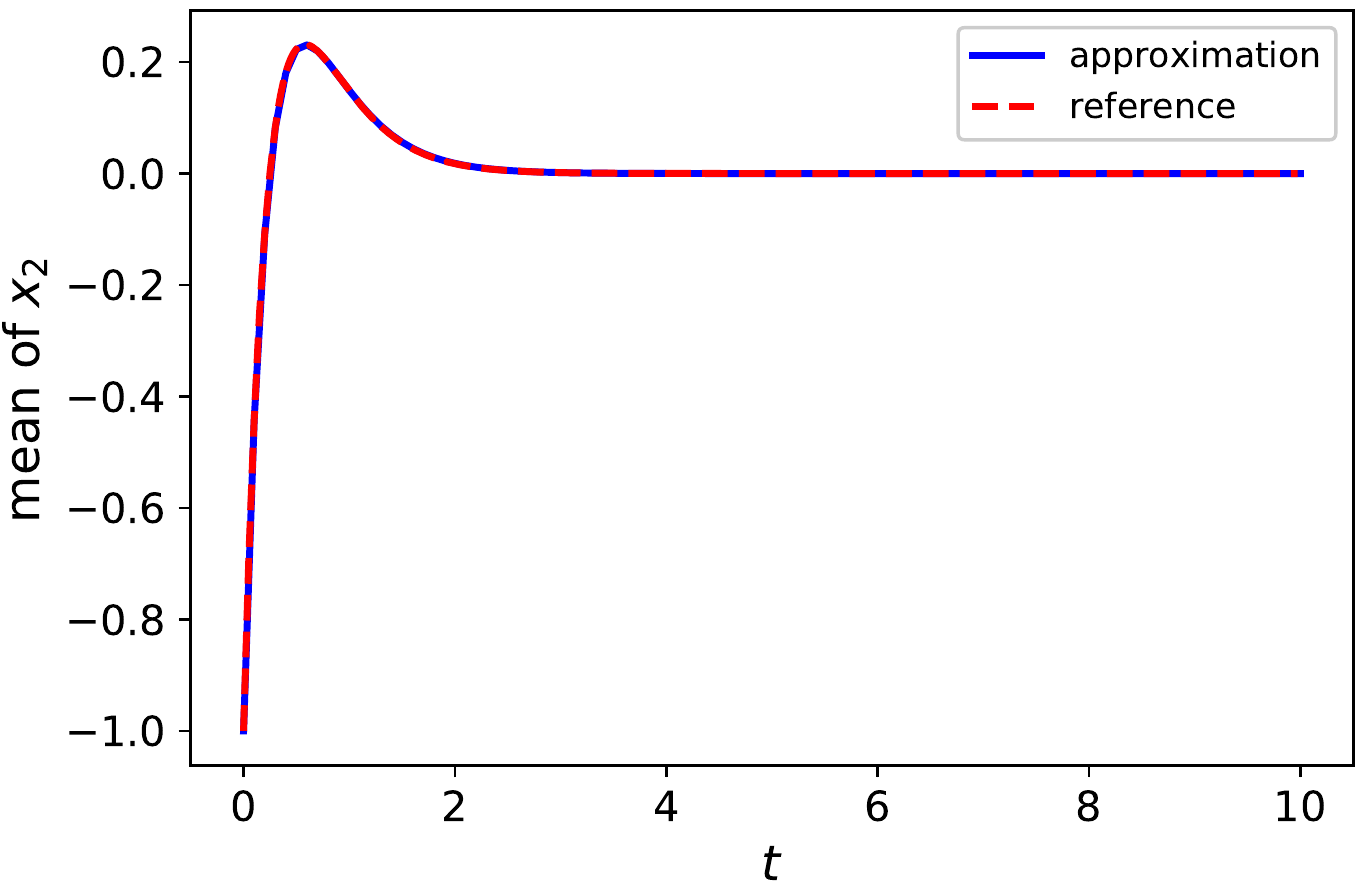}
			\caption{mean of $x_2$}
		\end{center}
	\end{subfigure}
	\begin{subfigure}[b]{0.48\textwidth}
		\centering
		\includegraphics[width=1.0\linewidth]{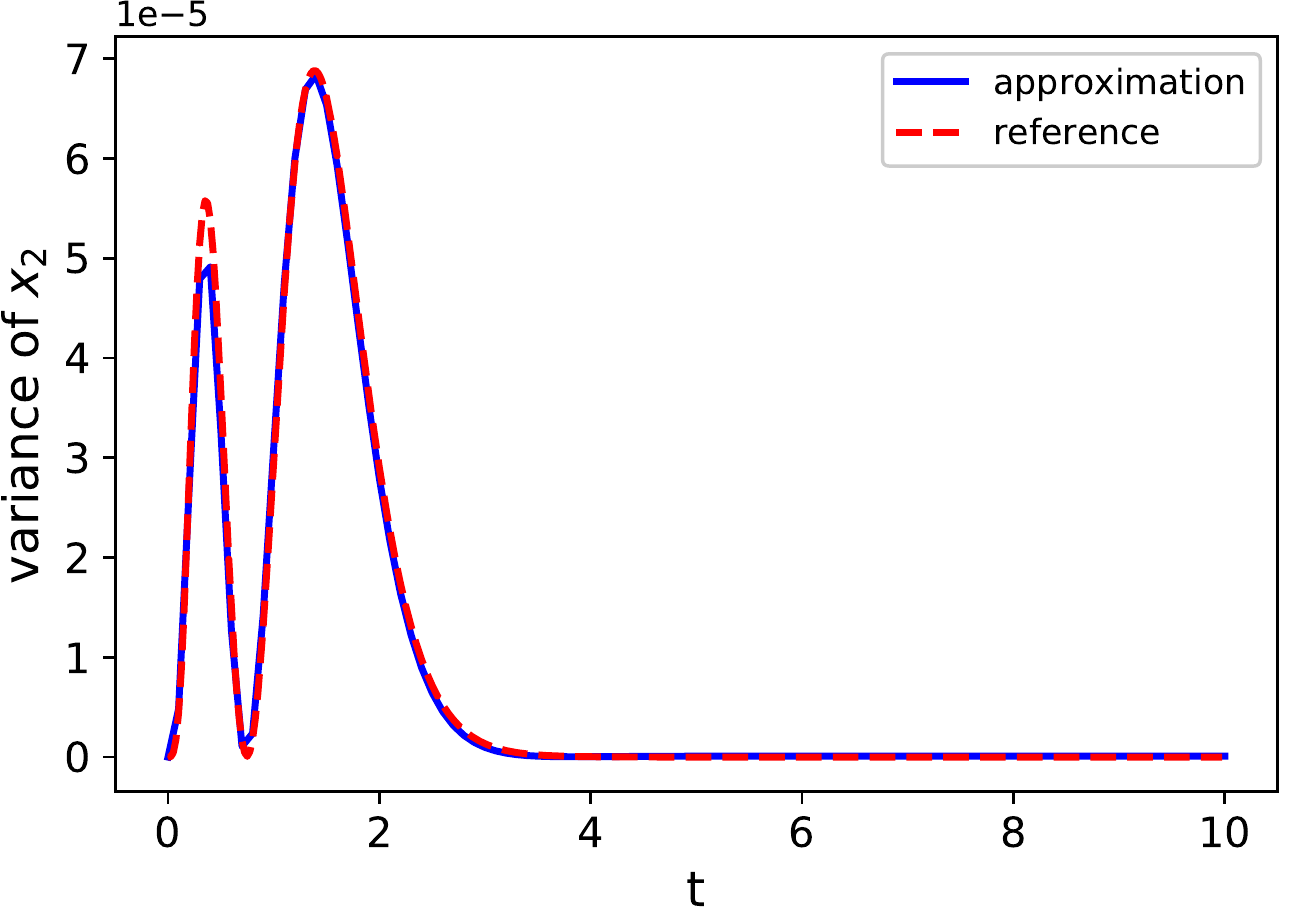}
		\caption{variance of $x_2$}
	\end{subfigure}	
	\caption{Mean and variance of the solution to Example 2 with $\x_0=(0, -1)$.}
	\label{fig:ex2}
\end{figure}

\begin{figure}[htb]
	\centering
	\begin{subfigure}[b]{0.48\textwidth}
		\begin{center}
			\includegraphics[width=1.0\linewidth]{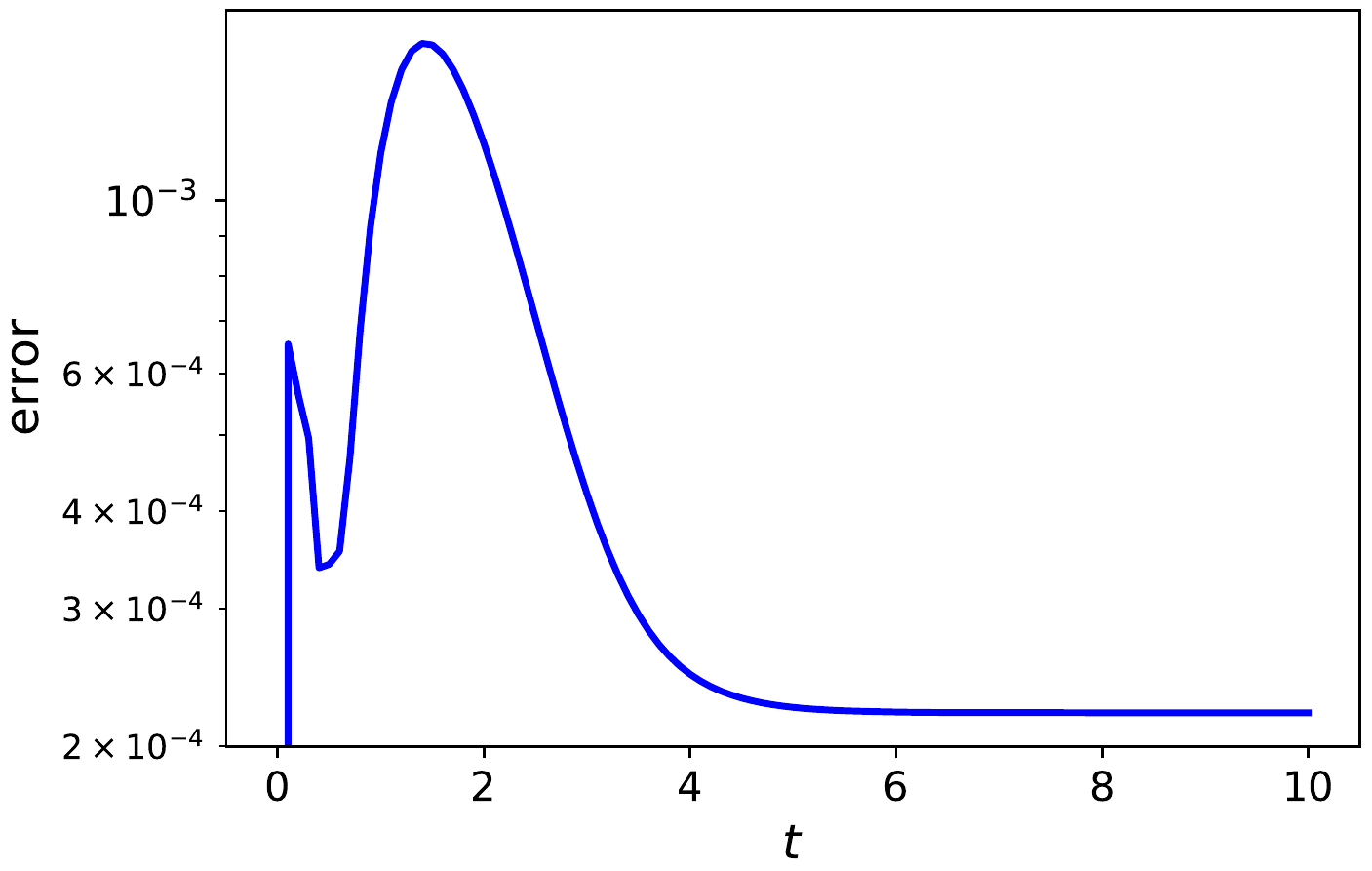}
			\caption{ error in the mean of $\x$ }
		\end{center}
	\end{subfigure}
	\begin{subfigure}[b]{0.48\textwidth}
		\centering
		\includegraphics[width=1.0\linewidth]{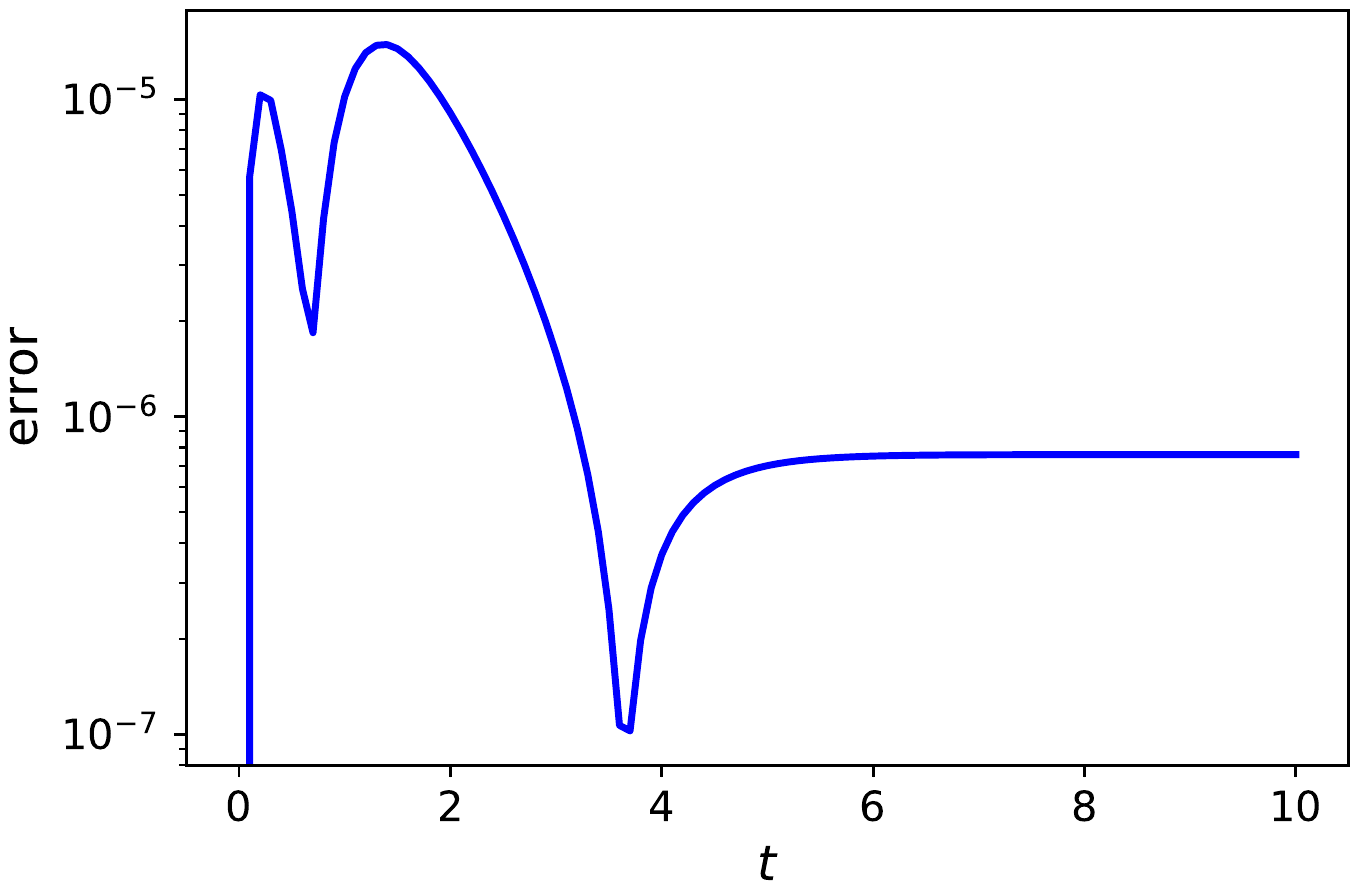}
		\caption{ error in the variance of $\x$ }
	\end{subfigure}
	\caption{The error of the mean and variance for Example 2.}
	\label{fig:ex2_mean_var_error}
\end{figure}

For UQ, we set the parameters $\aalpha$ to follow a multivariate uniform
distribution over $I_{\aalpha}$. In \figref{fig:ex2}, the approximated
mean and variance,  computed using a tensor product of five-point
Gauss-Legendre quadrature, are presented. Good match
between the DNN approximation and the reference is observed. In
\figref{fig:ex2_mean_var_error}, we see that the error in the mean and
the variance is around $10^{-3}$ and $10^{-5}$, respectively. 

\subsection{Example 3: Nonlinear Random Oscillation}

We now consider a nonlinear system of ODEs,
\begin{equation}
	\begin{split}
		\frac{dx_1}{dt} & = x_2, \\
		\frac{dx_2}{dt} & = -\alpha_1 x_2 - \alpha_2 \sin x_1,
	\end{split}
\end{equation}
where $\aalpha=(\alpha_1,\alpha_2)\in I_{\aalpha}=[0, 0.4] \times [8.8, 9.2]$ and $I_\x=[-\pi, \pi]\times [-2\pi, 2\pi]$.

\begin{figure}[htb]
	\centering
	\begin{subfigure}[b]{0.48\textwidth}
		\begin{center}
			\includegraphics[width=1.0\linewidth]{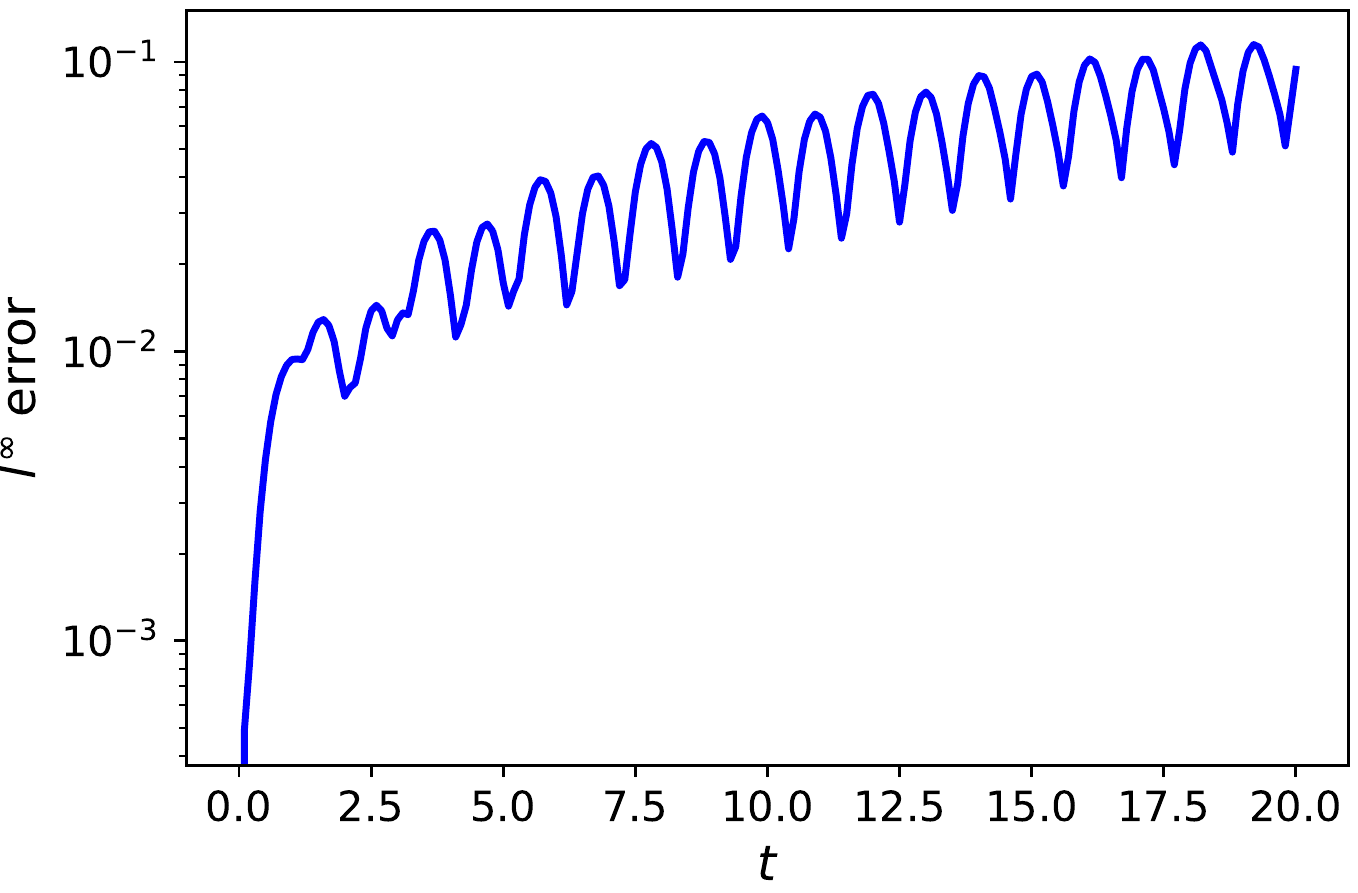}
			\caption{$l^\infty$ error}
		\end{center}
	\end{subfigure}
	\begin{subfigure}[b]{0.48\textwidth}
		\centering
		\includegraphics[width=1.0\linewidth]{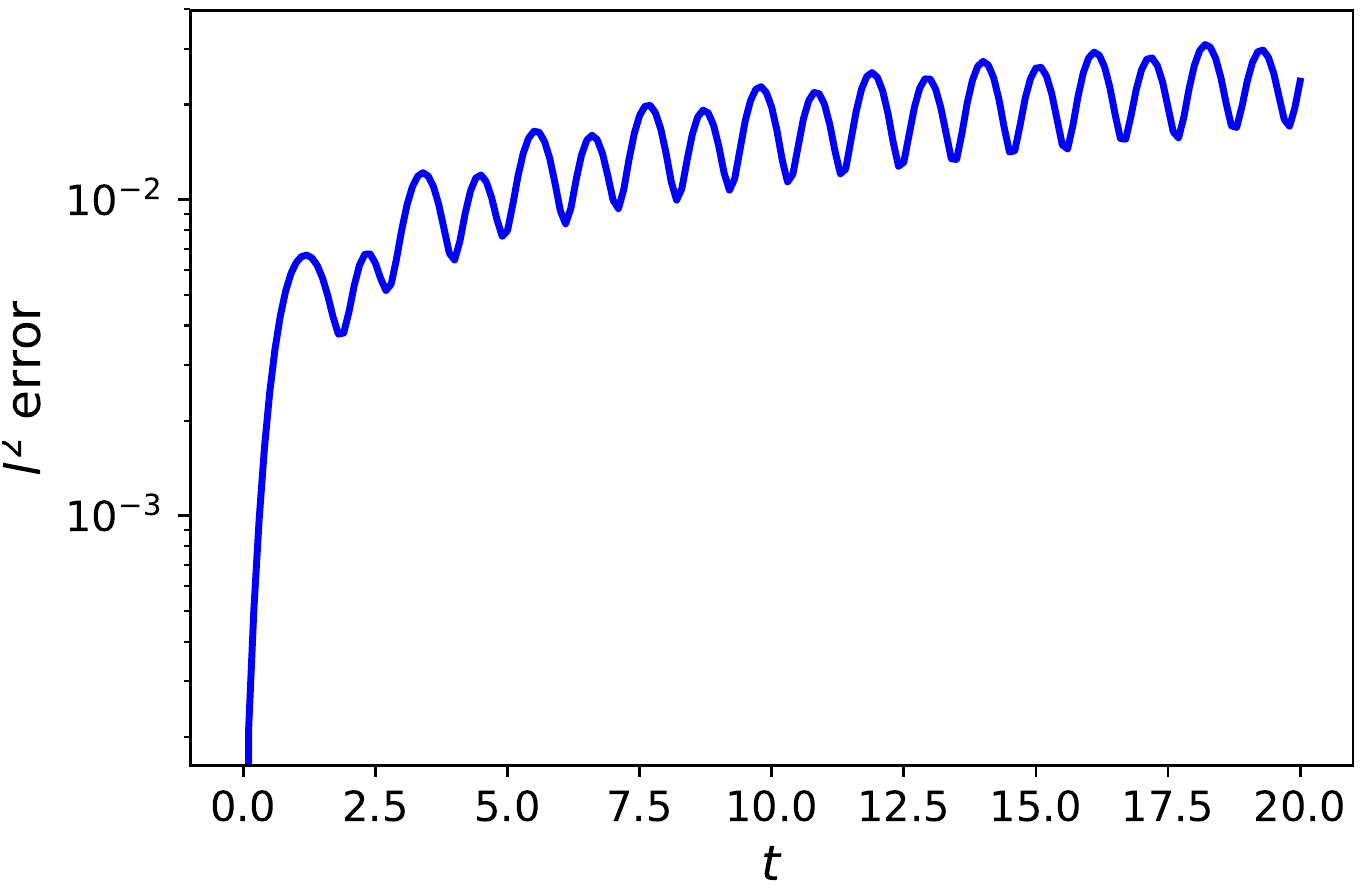}
		\caption{$l^2$ error}
	\end{subfigure}	
	\caption{The $l^\infty$ (left) and $l^2$ (right) errors of sample trajectories of Example 3 with $\x_0=(-1.193, -3.876)$.}
	\label{fig:ex3_error}
\end{figure}

The flow map is approximated by a DNN with structure $(3, 40)$. The
$l^\infty$ and $l^2$ errors for $1000$ sample trajectories with
$\x_0=(-1.193, -3.876)$ and $t=20$ are presented in
\figref{fig:ex3_error}. The error grows with time and stay around
$10^{-2}$. The oscillation is due the oscillatory behavior of the
solution. 

\begin{figure}[htb]
	\centering
	\begin{subfigure}[b]{0.48\textwidth}
		\begin{center}
			\includegraphics[width=1.0\linewidth]{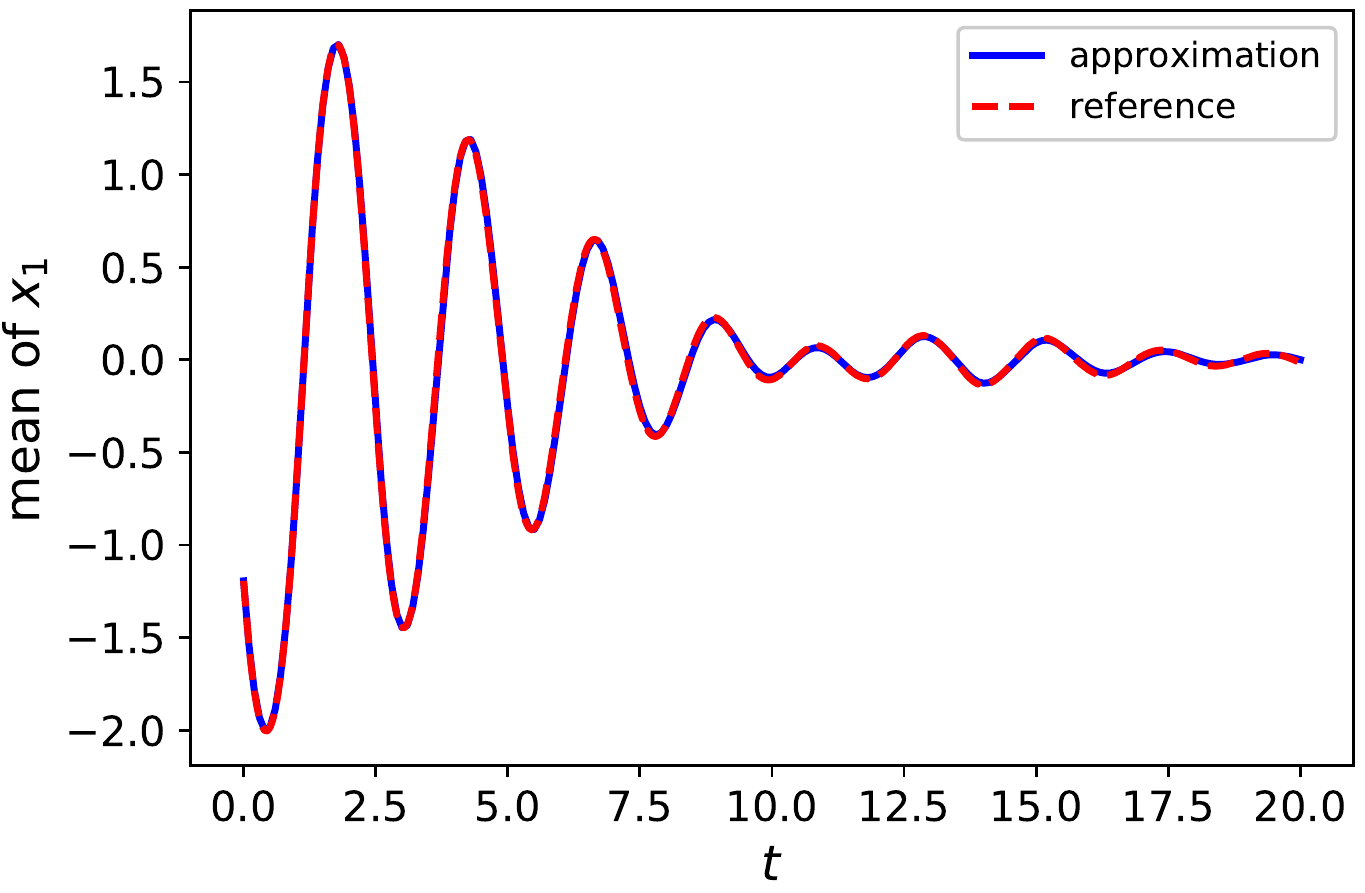}
			\caption{mean of $x_1$}
		\end{center}
	\end{subfigure}
	\begin{subfigure}[b]{0.48\textwidth}
		\centering
		\includegraphics[width=1.0\linewidth]{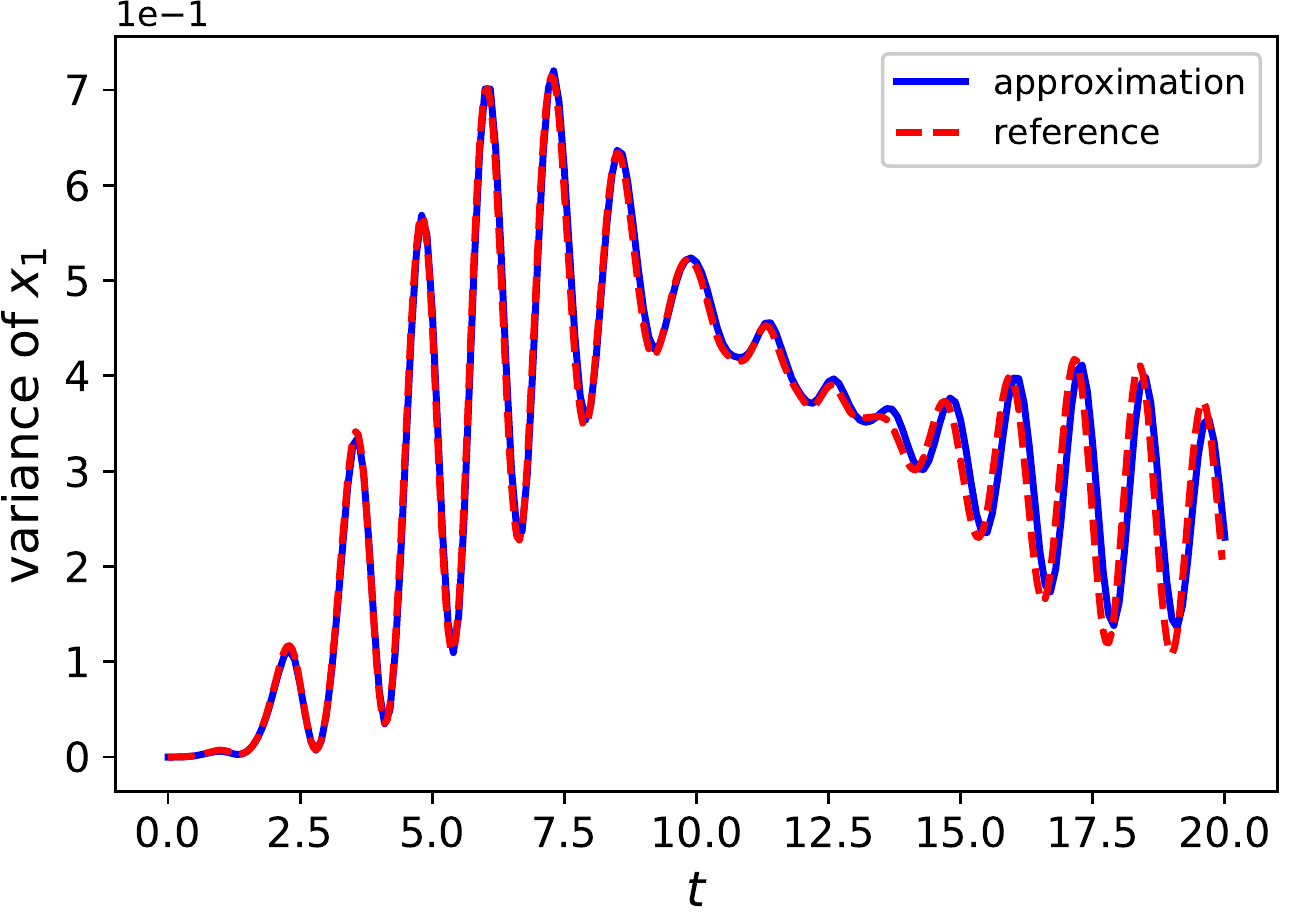}
		\caption{variance of $x_1$}
	\end{subfigure}
	\begin{subfigure}[b]{0.48\textwidth}
		\begin{center}
			\includegraphics[width=1.0\linewidth]{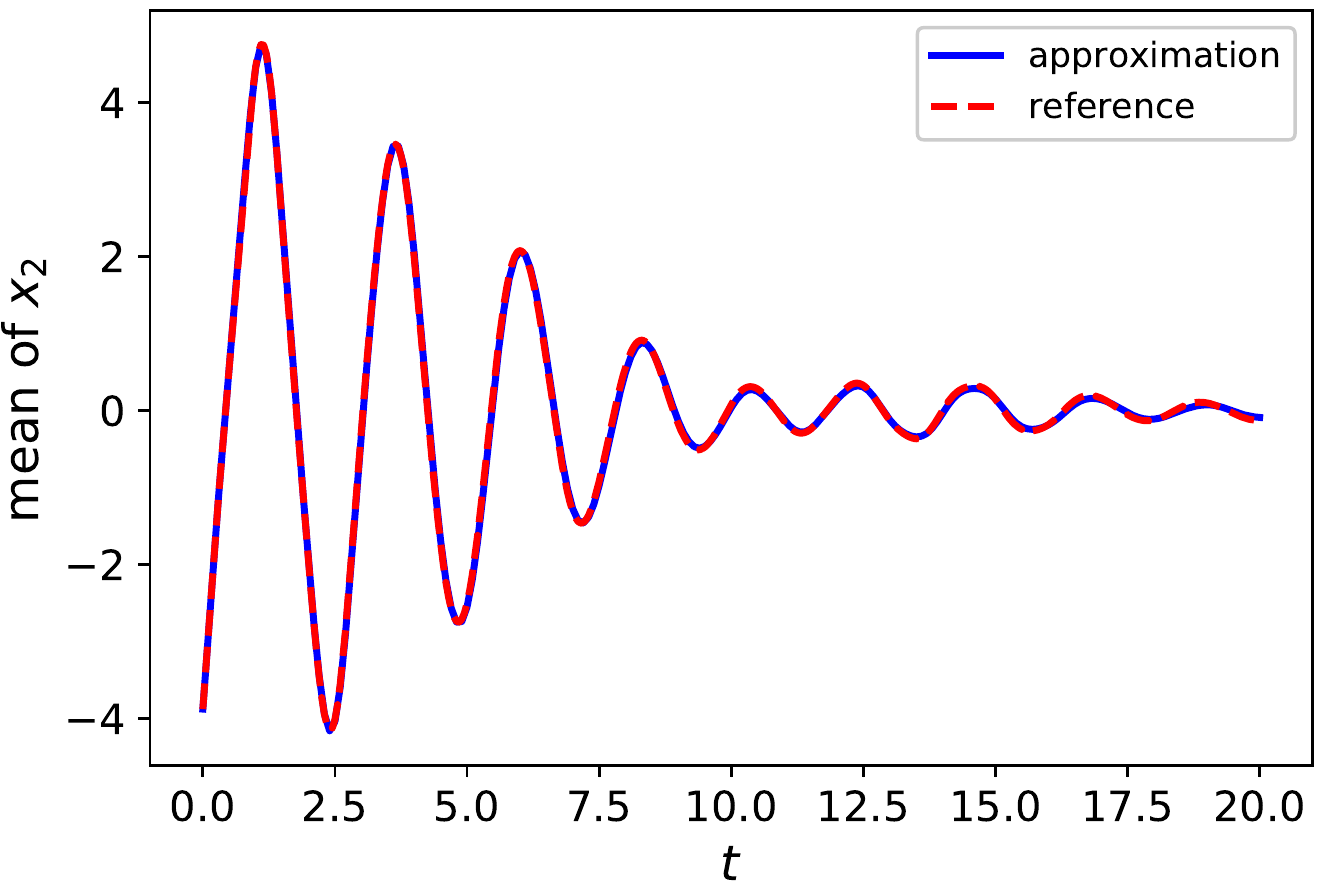}
			\caption{mean of $x_2$}
		\end{center}
	\end{subfigure}
	\begin{subfigure}[b]{0.48\textwidth}
		\centering
		\includegraphics[width=1.0\linewidth]{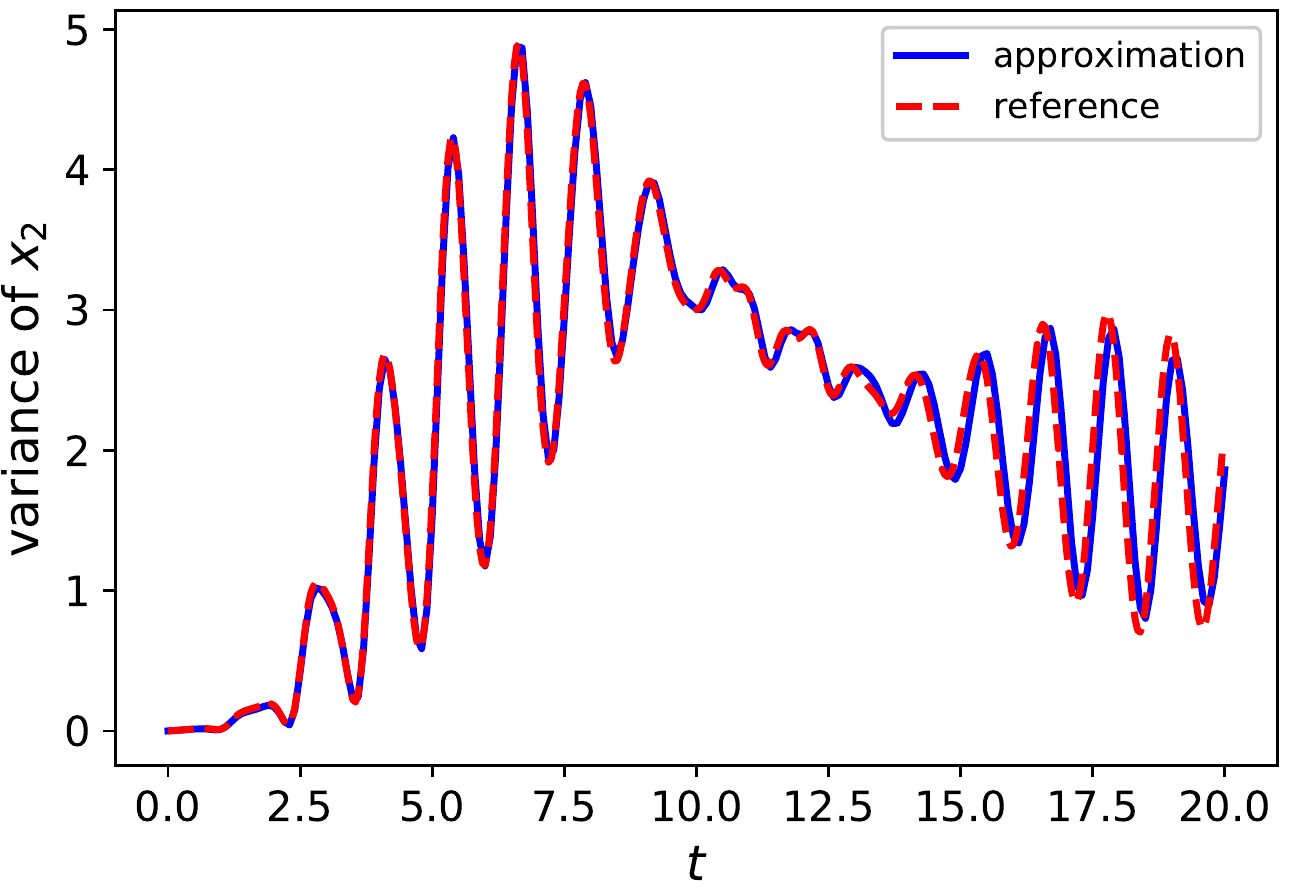}
		\caption{variance of $x_2$}
	\end{subfigure}
	
	\caption{Mean (left column) and variance (right column) of the solution to Example 3 with $\x_0=(-1.193, -3.876)$.}
	\label{fig:ex3}
\end{figure}

The mean and variance are computed in the same way as in Example 2. We present the results in \figref{fig:ex3}. A good performance is observed for the approximation for both the mean and the variance. For the variance approximation, a slight deviation is observed after $t=12.5$, i.e., $n=125$. This is due to the accumulation of the error as shown in \eqref{eq:estimate_var}. We further plot the propagation of the error of the variance approximation in \figref{fig:ex3_error}. The error grows slowly and oscillates with respect to $t$.

\begin{figure}[htb]
	\centering
	\begin{subfigure}[b]{0.48\textwidth}
		\begin{center}
			\includegraphics[width=1.0\linewidth]{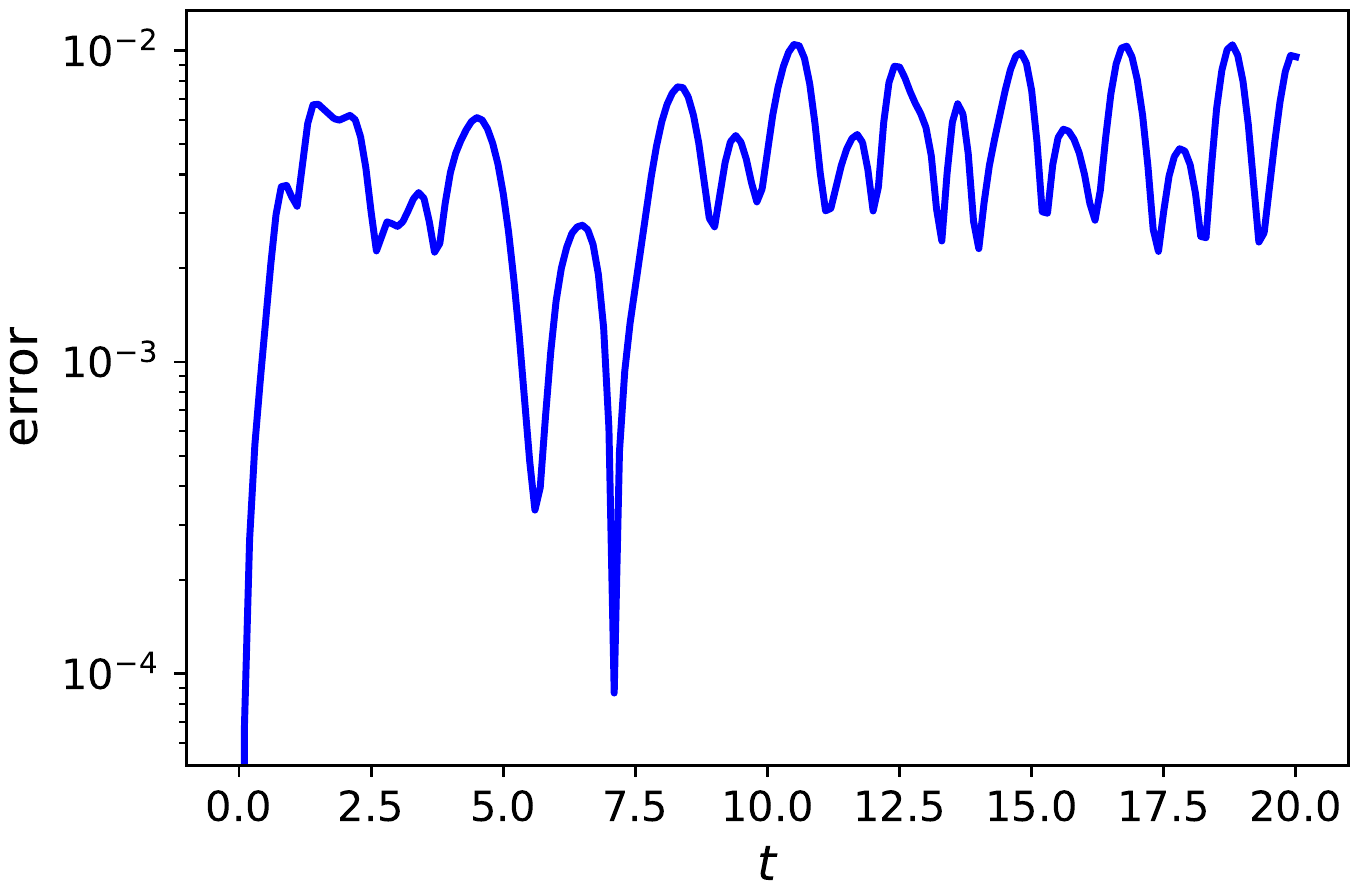}
			\caption{error of mean}
		\end{center}
	\end{subfigure}
	\begin{subfigure}[b]{0.48\textwidth}
		\centering
		\includegraphics[width=1.0\linewidth]{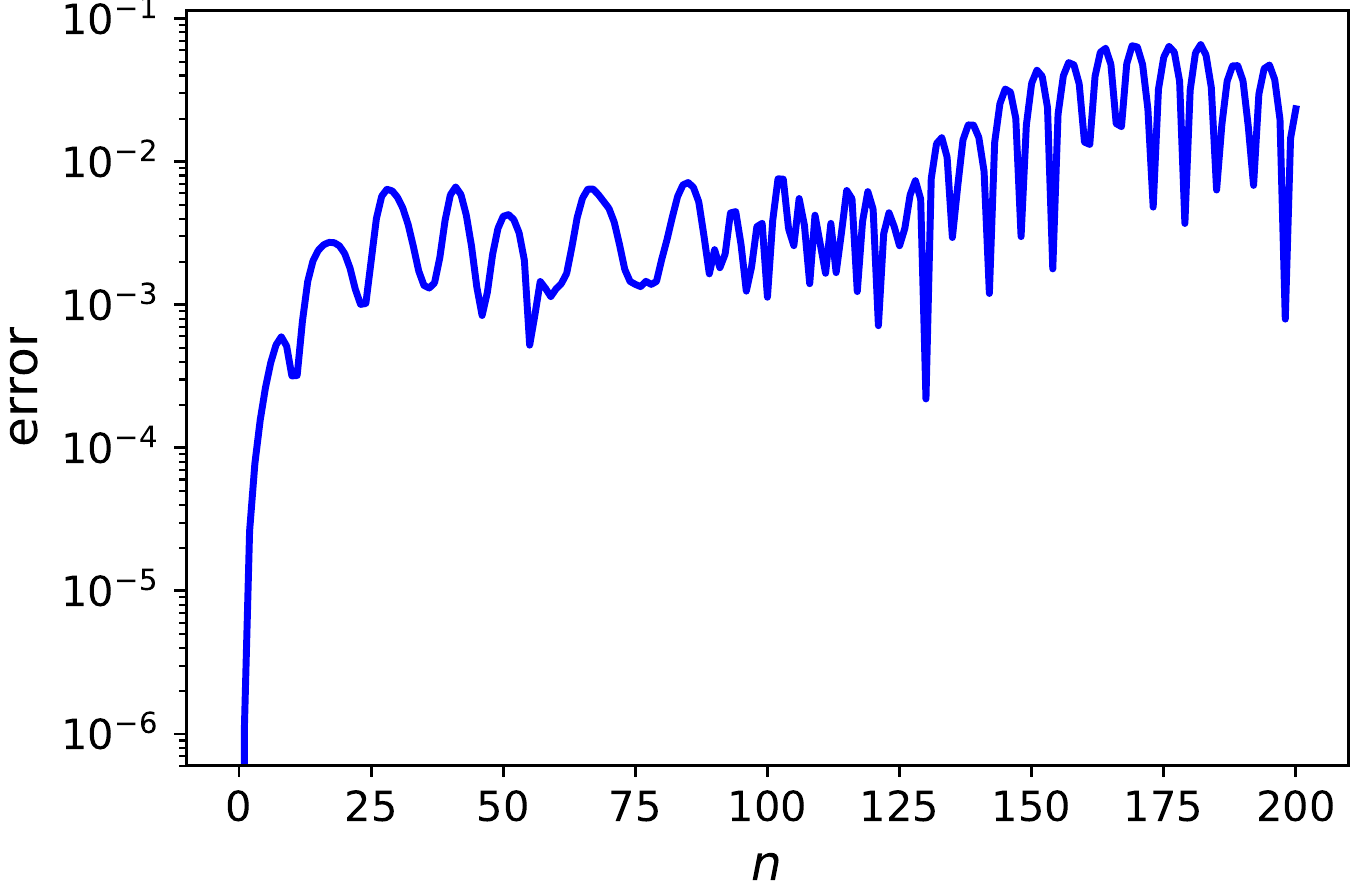}
		\caption{error of variance}
	\end{subfigure}	
	\caption{Propagation of errors in the mean and variance of the solution to Example 3 with $\x_0=(-1.193, -3.876)$.}
	\label{fig:ex3_error_mean_var}
\end{figure}

\subsection{Example 4: Cell Signaling Cascade}

The last example is nonlinear system with moderately high dimensional
parameter inputs. It is a mathematical model for autocrine
cell-signaling loop developed in \cite{shvartsman2002autocrine} in the
following form
\begin{equation}
\label{eq:CSC}
	\begin{split}
		\frac{d e_{1p}}{dt}&=\frac{I}{1+G e_{3p}} \frac{V_{\max, 1}(1-e_{1p})}{K_{m,1}+(1-e_{1p})}-\frac{V_{\max, 2} e_{1p}}{K_{m, 2}+e_{1p}},\\
		\frac{d e_{2p}}{dt}&=\frac{V_{\max, 3} e_{1p}(1-e_{2p})}{K_{m, 3} +(1-e_{2p})}-\frac{V_{\max, 4} e_{2p}}{K_{m, 4} + e_{2p}},\\
		\frac{d e_{3p}}{dt}&=\frac{V_{\max, 5} e_{2p}(1-e_{3p})}{K_{m, 5} +(1-e_{3p})}-\frac{V_{\max, 6} e_{3p}}{K_{m, 6} + e_{3p}}.
	\end{split}
\end{equation}
The state variables $e_{1p},\, e_{2p}$, and $e_{3p}$ denote the
dimensionless concentrations of the active form of the enzymes. This
model contains $13$ (random) parameters: $K_{m, 1-6}$, $V_{\max, 1-6}$, and
$G$, and a tuning parameter $I$ with range $[0, 1.5]$.
%
In \cite{shvartsman2002autocrine}, the parameters take the following
values $K_{m, 1-6}=0.2$, $V_{\max, 1}=0.5$, $V_{\max, 2}=0.15$,
$V_{\max, 3}=0.15$, $V_{\max, 4}=0.15$, $V_{\max, 5}=0.25$, $V_{\max,
  6}=0.05$, and $G=2$. Here, we use these values are the mean values
for the parameters and assume all parameters are independently and
uniformly distributed in a hypercube of $\pm 10\%$ around the mean values.

Each concentration should fall between $0$ and $1$ and hence we take
$I_\x=[0, 1]^3$. Moreover, to ensure the output of the DNN falls in
this physical bound, we add an activation function
$\sigma_{\text{output}}=\tanh(x)$ on each output node. The fully
connected block in the network employed here has a structure with $3$
layers and $200$ nodes each layer. For illustration purpose, we
calculate the mean and variance of the state variables with respect to
the random parameters $K_{m, 1}$, $K_{m,4}$, $V_{\max, 2}$, and
$V_{\max, 5}$. For other parameters, we assign the aforementioned
nominal values and treat them as deterministic. The tuning parameter
$I$ is taken to be $0.48$.  

\begin{figure}[htb]
	\centering
	\begin{subfigure}[b]{0.48\textwidth}
		\begin{center}
			\includegraphics[width=1.0\linewidth]{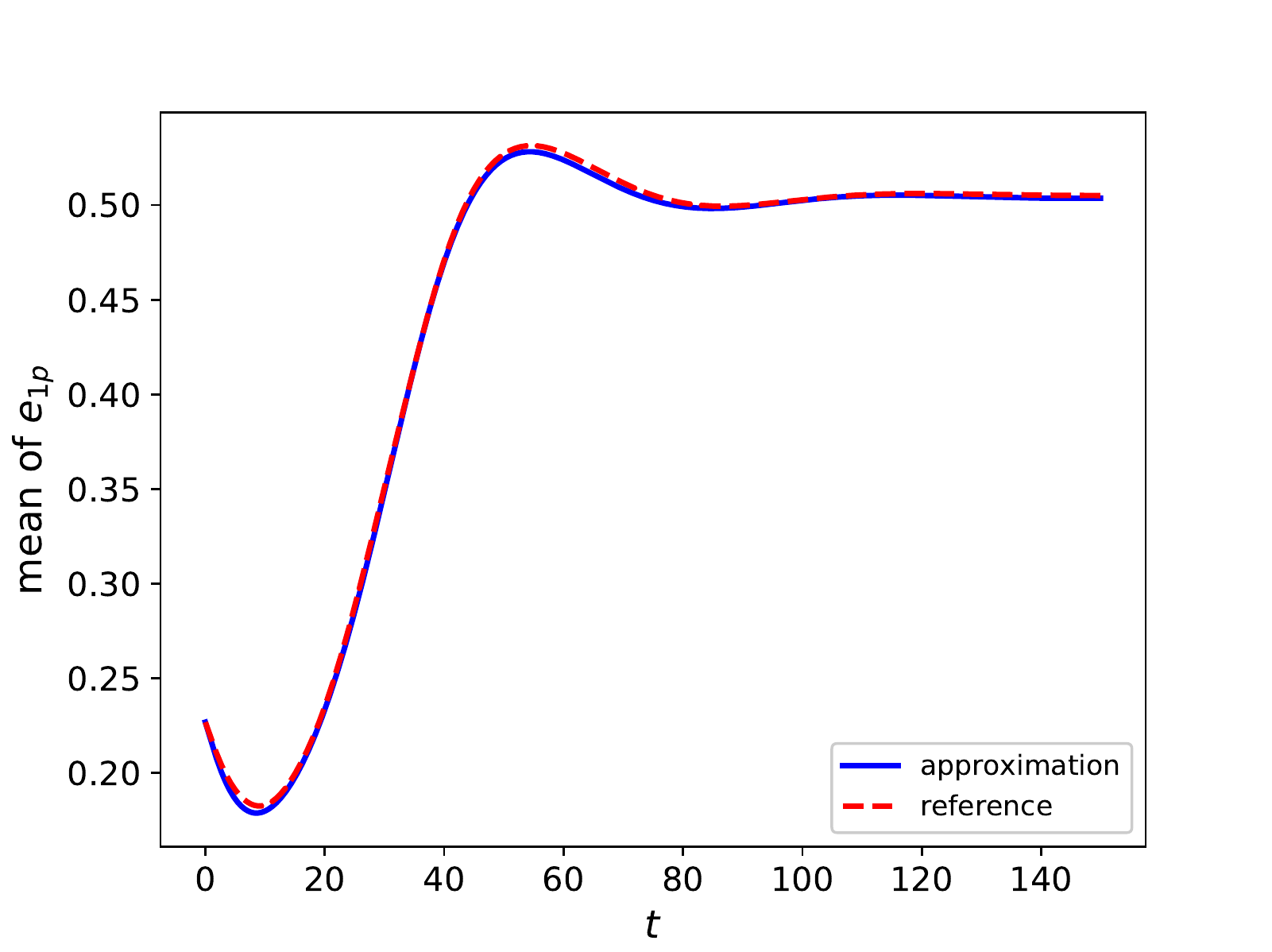}
			\caption{mean of $e_{1p}$}
		\end{center}
	\end{subfigure}
	\begin{subfigure}[b]{0.48\textwidth}
		\centering
		\includegraphics[width=1.0\linewidth]{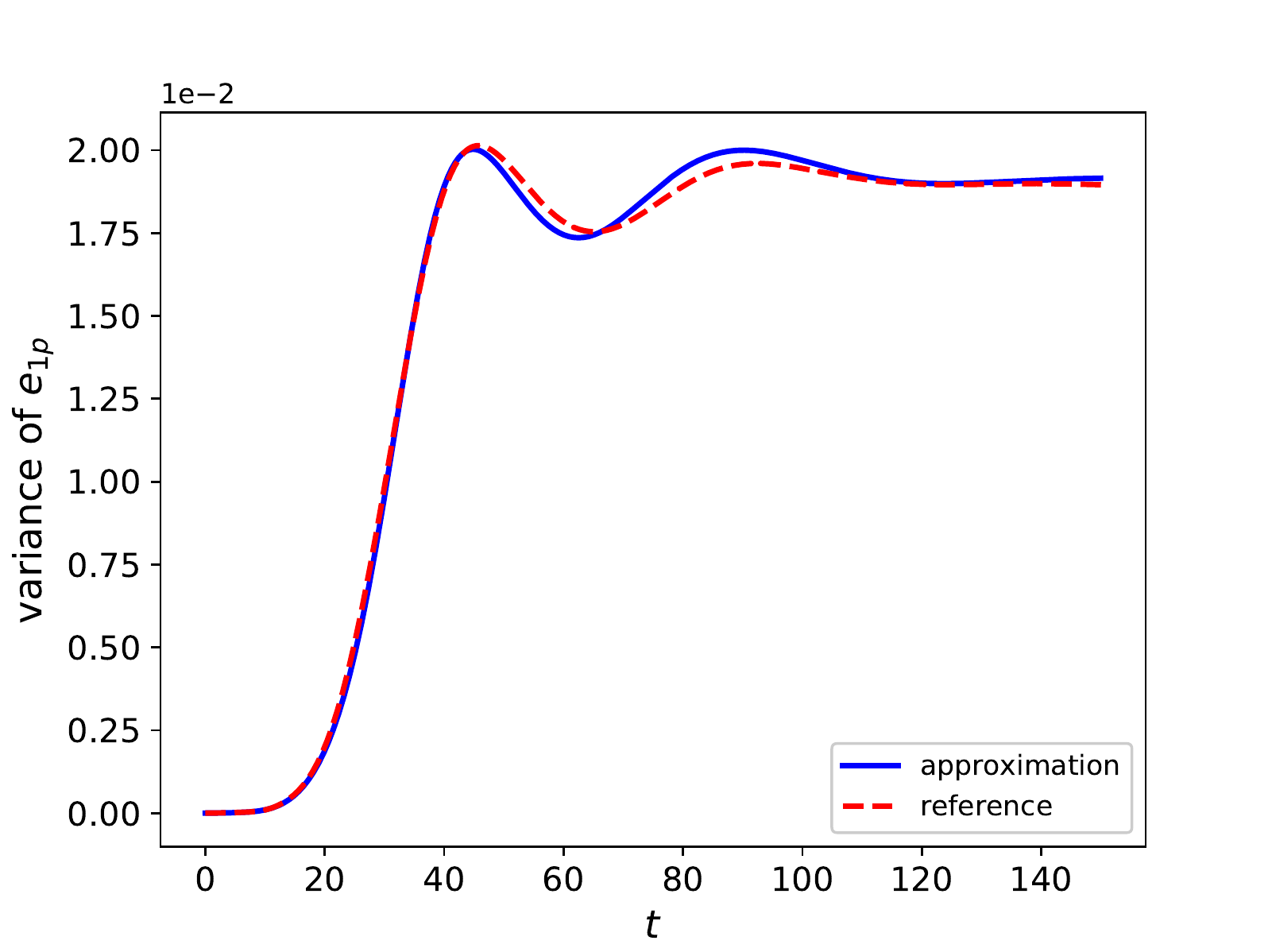}
		\caption{variance of $e_{1p}$}
	\end{subfigure}
	
	\begin{subfigure}[b]{0.48\textwidth}
		\begin{center}
			\includegraphics[width=1.0\linewidth]{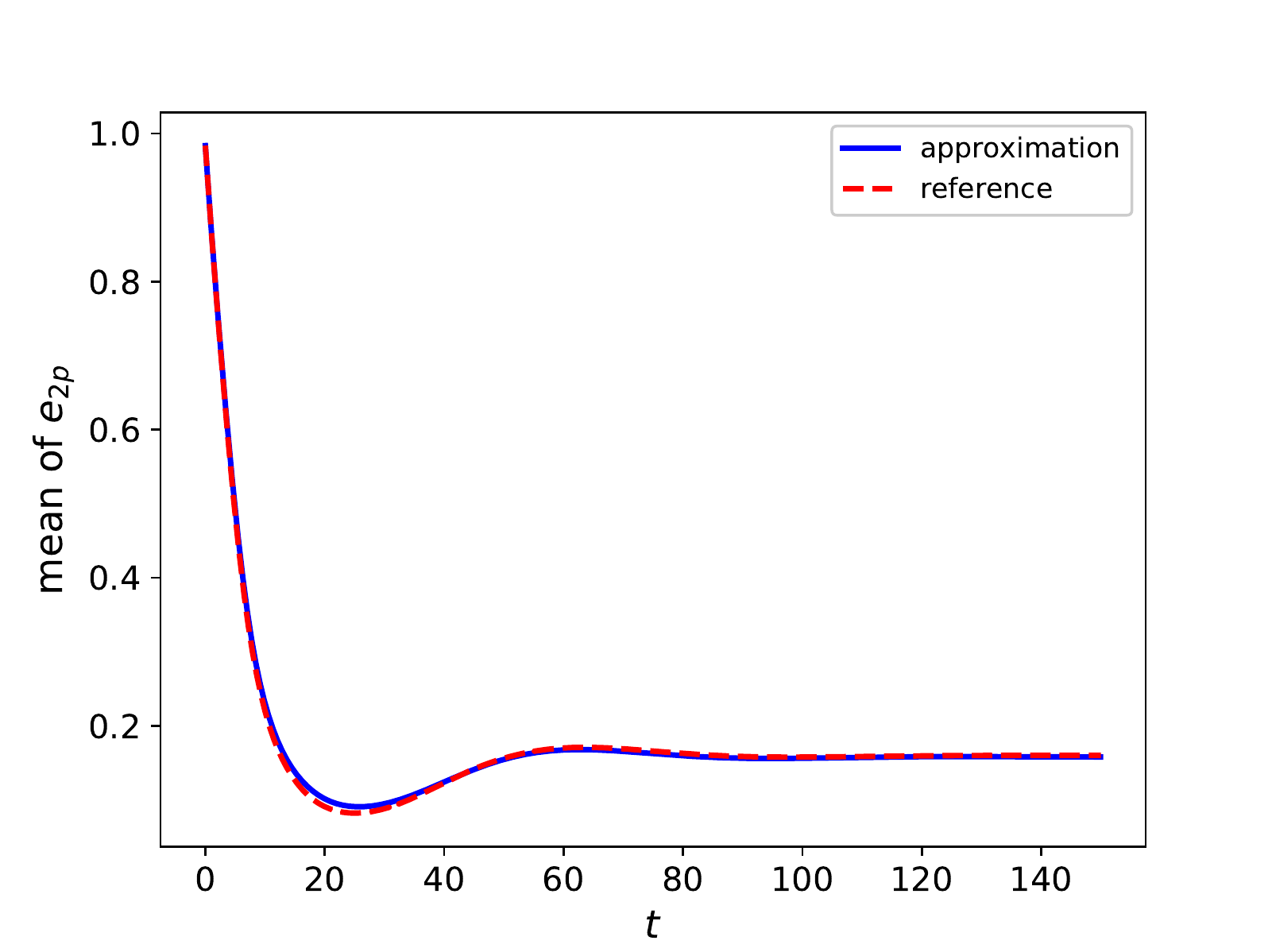}
			\caption{mean of $e_{2p}$}
		\end{center}
	\end{subfigure}
	\begin{subfigure}[b]{0.48\textwidth}
		\centering
		\includegraphics[width=1.0\linewidth]{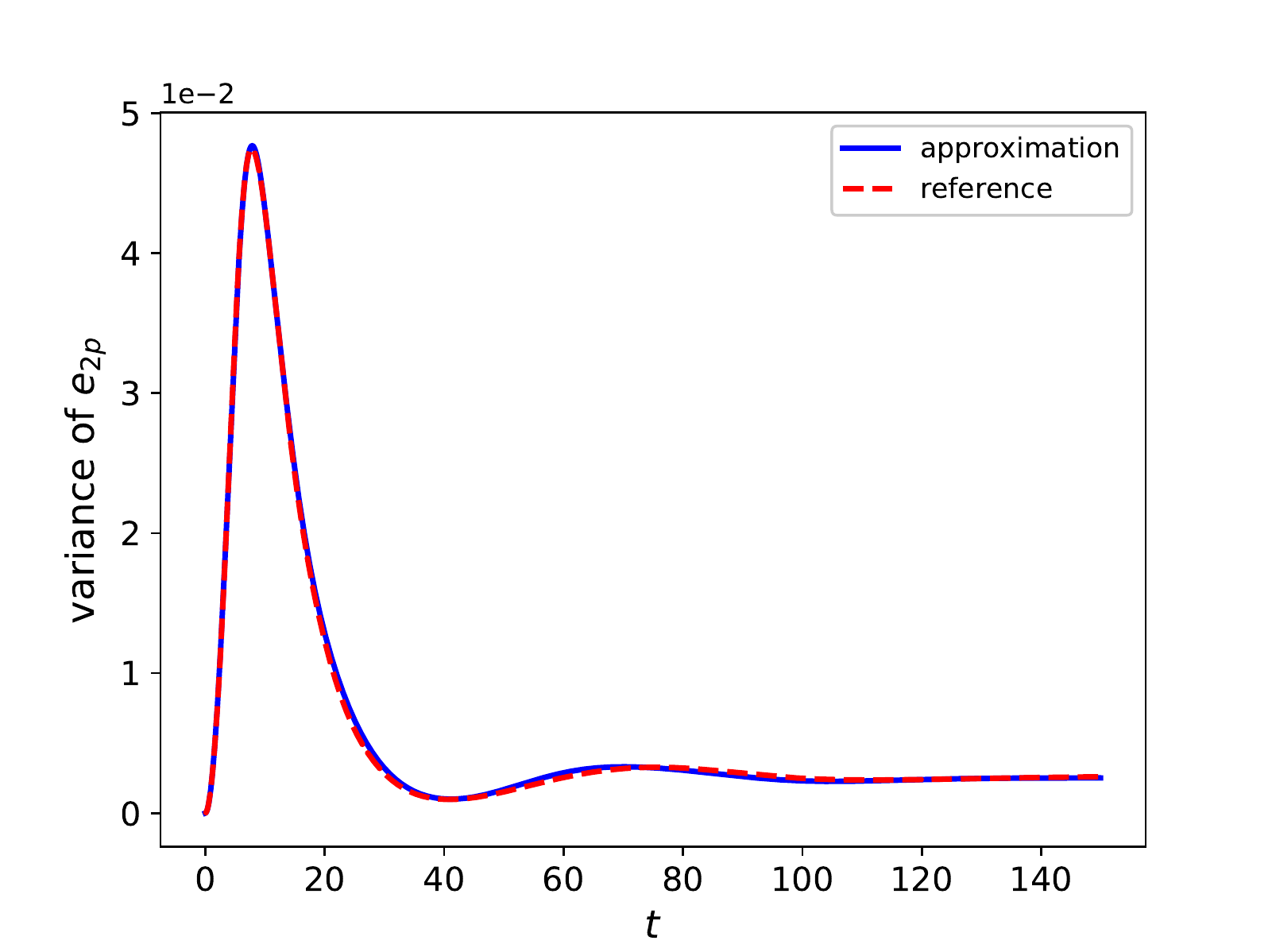}
		\caption{variance of $e_{2p}$}
	\end{subfigure}
	
	\begin{subfigure}[b]{0.48\textwidth}
		\begin{center}
			\includegraphics[width=1.0\linewidth]{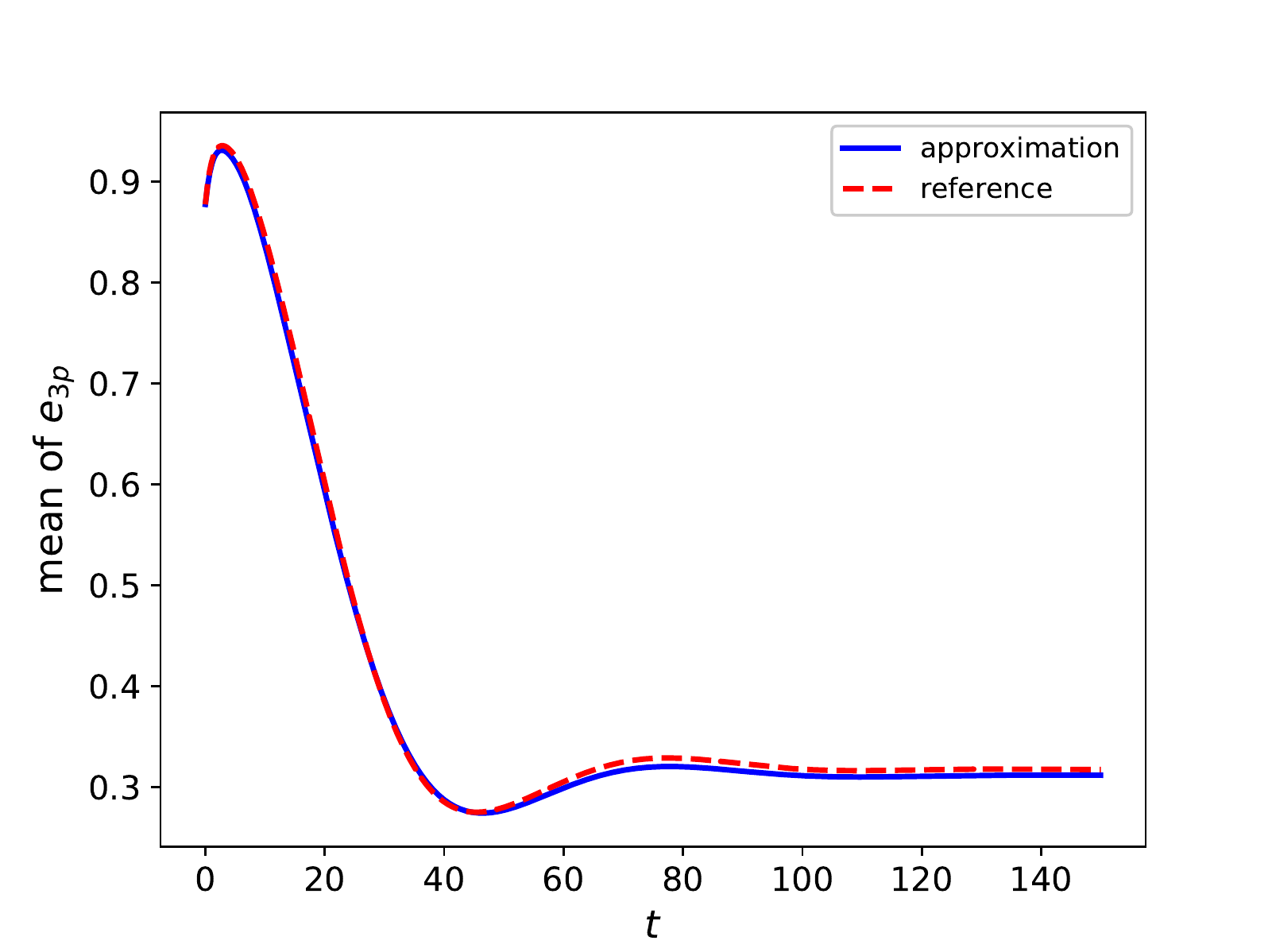}
			\caption{mean of $e_{3p}$}
		\end{center}
	\end{subfigure}
	\begin{subfigure}[b]{0.48\textwidth}
		\centering
		\includegraphics[width=1.0\linewidth]{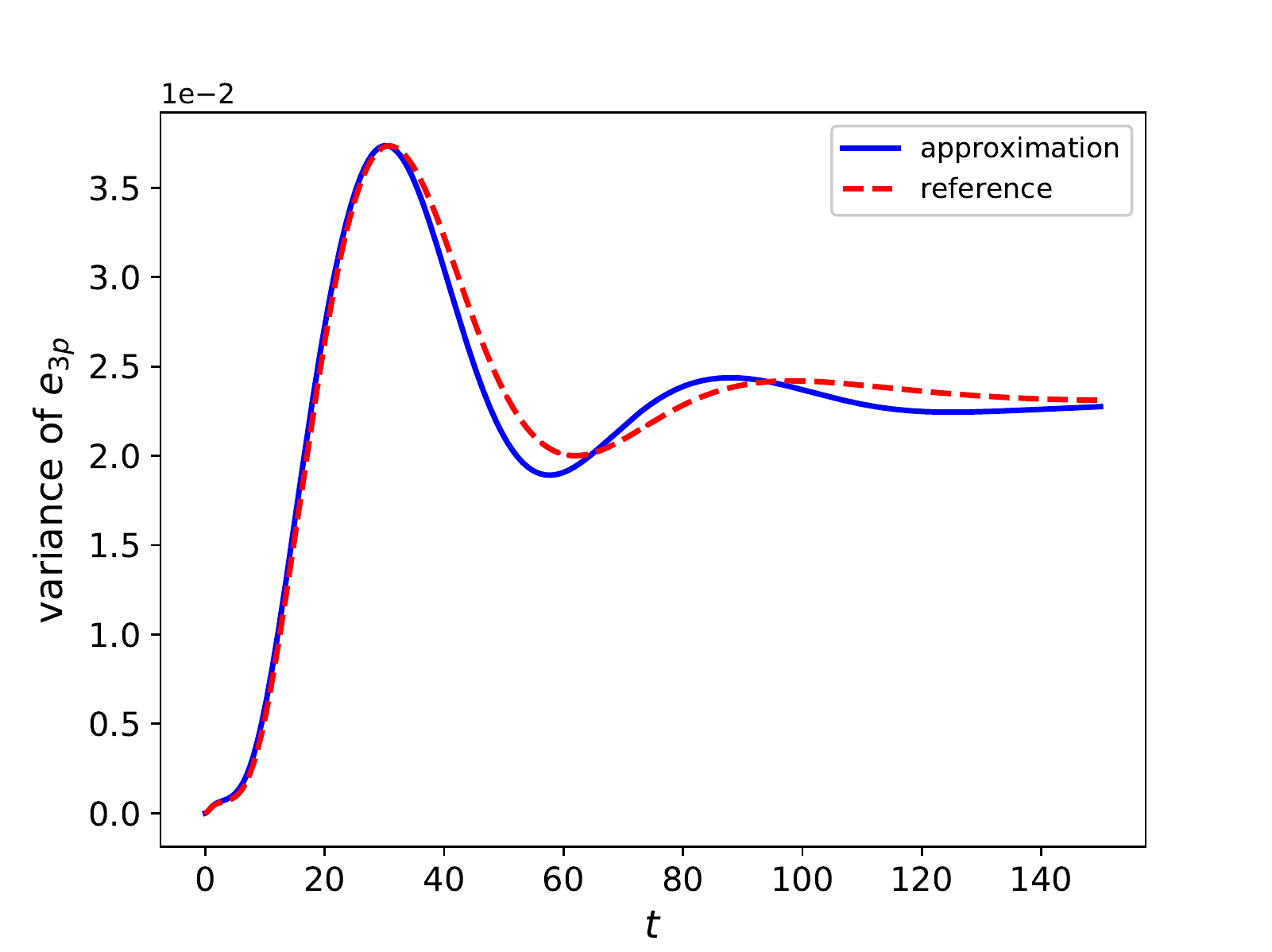}
		\caption{variance of $e_{3p}$}
	\end{subfigure}
	\caption{Mean (left column) and variance (right column) of the solution to Example 4.}
	\label{fig:ex4}
\end{figure}

After the training is finished, we march forward for $n=1,400$ steps
with the initial condition $\x_0=(0.22685145, 0.98369158, 0.87752945)$
and compute the mean and variance with a tensor product of five-point
Gaussian quadrature. In \figref{fig:ex4}, we present the approximated
mean and variance. Given such long-time simulation, the approximation
agrees with the reference. 

\begin{figure}[htb]
	\centering
	\includegraphics[width=1.0\linewidth]{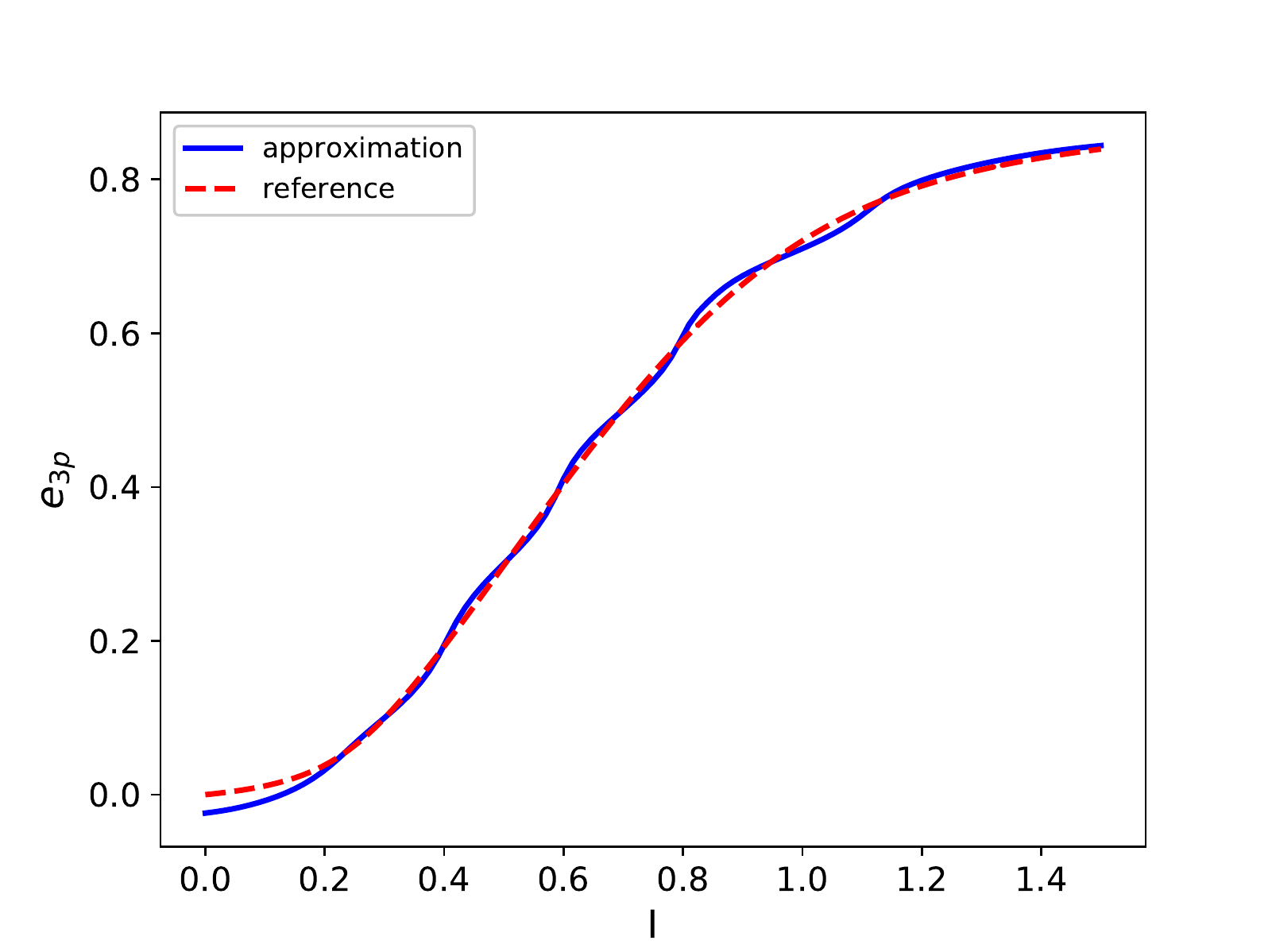}	
	\caption{Response curve of the steady state of $e_{3p}$ with respect to $I$ for Example 4.}
	\label{fig:ex4_e3p_I}
\end{figure}

For this example, the response curve of $e_{3p}$ with respect to the tuning parameter $I$ is of particular interest in practice. We examine such curve at the steady state of $e_{3p}$. To this end, we fix all the other parameters at their mean value and let $I$ vary in $[0, 1.5]$. To reach the steady state, without solving the true governing equations for long time, we march the DNN model forward for $2000$ steps as in \eqref{prediction}. \figref{fig:ex4_e3p_I} presents the resulting response curve. There is a good agreement between the approximation and the reference response curve observed.

\section{Conclusion} \label{sec:conclusions}
We presented a numerical framework for discovering unknown parameterized dynamical systems, using observational data and deep neural networks (DNN). The network structure is a modification to the residual network, which was shown to be effective to discover unknown deterministic dynamical systems in recent work \cite{QinWuXiu2019}.
Our method allows one to not only create accurate neural network model to approximate the unknown dynamics but also conduct efficient uncertainty quantification.

\section{Acknowledgments}
Sandia National Laboratories is a multi-mission laboratory managed and operated by National Technology and Engineering Solutions of Sandia, LLC., a wholly owned subsidiary of Honeywell International, Inc., for the U.S. Department of Energy's National Nuclear Security Administration under contract DE-NA-0003525. The views expressed in the article do not necessarily represent the views of the U.S. Department of Energy or the United States Government.

\clearpage
\bibliographystyle{siamplain}
\bibliography{neural,LearningEqs,UQ}

%
%

\end{document}